\documentclass[aap]{imsart}

\RequirePackage{amsthm,amsmath,amsfonts,amssymb}
\RequirePackage[authoryear]{natbib}%% uncomment this for author-year citations

\startlocaldefs
\numberwithin{equation}{section}
\theoremstyle{plain}
\newtheorem{theorem}{Theorem}[section]
\newtheorem{lemma}[theorem]{Lemma}
\newtheorem{proposition}[theorem]{Proposition}

\theoremstyle{remark}
\newtheorem{remark}[theorem]{Remark}

\newtheorem{definition}[theorem]{Definition}
\newtheorem{condition}[theorem]{Condition}
\def \non {{\nonumber}}

\def \tilde {\widetilde}
\def \E {\mathbb{E}}
\def \P {\mathbb{P}}
\def \R {\mathbb{R}}
\def \N {\mathbb{N}}
\def \1 {\mathbf{1}}

\endlocaldefs

\begin{document}

\begin{frontmatter}
\title{A reverse ergodic theorem for inhomogeneous killed Markov chains and application to a new uniqueness result for reflecting diffusions}
%\title{A sample article title with some additional note\thanksref{T1}}
\runtitle{An ergodic theorem and uniqueness for reflecting diffusions}
%\thankstext{T1}{A sample of additional note to the title.}

\begin{aug}
%%%%%%%%%%%%%%%%%%%%%%%%%%%%%%%%%%%%%%%%%%%%%%%
%% Only one address is permitted per author. %%
%% Only division, organization and e-mail is %%
%% included in the address.                  %%
%% Additional information can be included in %%
\author[A]{\fnms{Cristina}~\snm{Costantini}\ead[label=e1]{c.costantini@unich.it}}
\and
\author[B]{\fnms{Thomas G.}~\snm{Kurtz}\ead[label=e2]{kurtz@math.wisc.edu}}
%%%%%%%%%%%%%%%%%%%%%%%%%%%%%%%%%%%%%%%%%%%%%%
%% Addresses                                %%
%%%%%%%%%%%%%%%%%%%%%%%%%%%%%%%%%%%%%%%%%%%%%%
\address[A]{Department of Economic Studies and INdAM Local Unit, University of Chieti-Pescara\printead[presep={,\ }]{e1}}

\address[B]{Department of Mathematics and Department of Statistics, University of Wisconsin-Madison\printead[presep={,\ }]{e2}}

\end{aug}

\begin{abstract}
Bass and Pardoux (1987) deduce from the Krein-Rutman 
theorem a reverse ergodic theorem for a sub-probability 
transition function, which turns out to be a key tool in proving  
uniqueness of reflecting Brownian motion in cones in Kwon and 
Williams (1991) and Taylor and Williams (1993). 
By a different approach, we are able to prove an analogous reverse ergodic theorem for 
a family of inhomogeneous sub-probability transition 
functions. 

This allows us to prove existence and uniqueness for a 
semimartingale diffusion process with varying, 
oblique direction of reflection, in a domain with one singular point 
that can be approximated, near the singular point, by a 
smooth cone, under natural, easily verifiable geometric 
conditions. 

Along the way we also show that if the 
reflecting Brownian motion in a smooth cone is a 
semimartingale then the 
parameter $\alpha$ of Kwon and Williams (1991) is strictly less 
than 1, thus partially extending 
the results of Williams (1985) to higher dimension. 
\end{abstract}



\end{frontmatter}
%%%%%%%%%%%%%%%%%%%%%%%%%%%%%%%%%%%%%%%%%%%%%%
%% Please use \tableofcontents for articles %%
%% with 50 pages and more                   %%
%%%%%%%%%%%%%%%%%%%%%%%%%%%%%%%%%%%%%%%%%%%%%%
%\tableofcontents

%%%%%%%%%%%%%%%%%%%%%%%%%%%%%%%%%%%%%%%%%%%%%%
%%%% Main text entry area:

\setcounter{equation}{0}

\section{Introduction}\label{sectionintro}

Let $E$ be a compact metric space, and let $Q(x,dy)$ be a 
sub-probability transition 
function on $E$, that is $Q(x,dy)$ satisfies all conditions for 
a probability transition function except that 
$Q(x,E)\leq 1$. As it is well known, $Q(x,dy)$ is the transition 
function of a killed Markov chain. In the proof of their Theorem 5.4,  
\cite{BP87} show that, if $Q$ satisfies the conditions of 
the Krein-Rutman theorem (Theorems 6.1 and 6.3 of 
\cite{KR50}), then, for any pair of continuous functions 
$f,g$, $g>0$, 
and any sequence of probability measures on $E$, $\{\mu_k\}$, 
\begin{equation}\lim_{k\rightarrow\infty}\frac {\int Q^kf(x)\,\mu_
k(dx)}{\int Q^kg(x)\,\mu_k(dx)}=C(f,g),\label{eq:BP}\end{equation}
where the constant $C(f,g)$ is independent of the sequence 
$\{\mu_k\}$. In particular, with $g\equiv 1$, 
\[\lim_{k\rightarrow\infty}\frac {\int Q^kf(x)\,\mu_k(dx)}{\int Q^
k\mathbf{1}(x)\,\mu_k(dx)}=C(f).\]
$(\ref{eq:BP})$ can be viewed as a {\em reverse ergodic }
{\em theorem for killed Markov chains}. Note that, tipically, 
both the numerator and the denominator in $(\ref{eq:BP})$ 
tend to zero. 
The result of \cite{BP87} is a key element in the proof of 
uniqueness of reflecting Brownian motion in a smooth cone, with radially 
constant direction of reflection, by \cite{KW91}, and in the proof of 
uniqueness of reflecting Brownian motion in a polyhedral cone, 
with constant directions of reflection on 
each face, by \cite{TW93} and \cite{DW96}.

Our first goal here is to extend the \cite{BP87} result 
to a sequence of 
compact metric spaces $E_0,E_1,E_2,\ldots$ and a sequence of 
sub-transition functions $Q_1,Q_2,\ldots$, with $Q_k$ governing 
transitions from $E_k$ to $E_{k-1}$, 
and give conditions under which 
\begin{equation}\lim_{k\rightarrow\infty}\frac {\int Q_kQ_{k-1}\cdots 
Q_1f(x)\mu_k(dx)}{\int Q_kQ_{k-1}\cdots Q_1g(x)\mu_k(dx)}=C(f,g),\label{eq:constlim}\end{equation}
where $C(f,g)$ is independent of $\{\mu_k\}$. We call 
$(\ref{eq:constlim})$ a {\em reverse ergodic theorem}
{\em for inhomogeneous killed Markov chains}. 
Note that, even in the case when $E_k=E$ for all $k$ and 
$Q_k$ converges, as $k$ goes to infinity, to 
a sub-transition function $Q$ on $E$,   
it is not in general possible 
to obtain the limit in $(\ref{eq:constlim})$ from the 
Krein-Rutman theorem. In fact this would essentially 
require exchanging the limits 
\[\lim_{l\rightarrow\infty}\lim_{k\rightarrow\infty}\frac {\int Q_{
l+k}Q_{l+k-1}\cdots Q_{l+1}f(x)\mu_k(dx)}{\int Q_{l+k}Q_{l+k-1}\cdots 
Q_{l+1}g(x)\mu_k(dx)}\]
and 
\[\lim_{k\rightarrow\infty}\lim_{l\rightarrow
\infty}\frac {\int Q_{l+k}Q_{l+k-1}\cdots Q_{l+1}f(x)\mu_k(dx)}{\int 
Q_{l+k}Q_{l+k-1}\cdots Q_{l+1}g(x)\mu_k(dx)}.\]
Rather than trying to reinforce the conditions of the 
Krein-Rutman theorem, we provide new conditions under 
which $(\ref{eq:constlim})$ holds (Theorem 
\ref{th:ergodic}). Our conditions  
are uniform lower bounds which have a clear 
probabilistic meaning and can be verified in many 
contexts. 

In the second part of this paper we use our 
reverse ergodic theorem for inhomogeneous killed Markov 
chains to prove uniqueness for a semimartingale reflecting 
diffusion process with varying, oblique direction of reflection, 
in a curved domain with only one singular point, that, in 
a neighborhood of the singular point, can be approximated 
by a smooth cone.  
Although one expects that such a diffusion, for a 
short time after it leaves the singular point, can be 
approximated by a 
reflecting Brownian motion in the cone (we actually prove this 
in a rigorous sense: see Lemma \ref{th:scaling}), it is not 
clear how the uniqueness result we seek might follow from the 
\cite{KW91} result, essentially due to a limit exchange 
problem like the one mentioned above. 
We prove existence and uniqueness of a semimartingale 
reflecting diffusion process as above under very 
natural, geometric conditions (Conditions \ref{D}, \ref{G} 
and \ref{A}): Besides mild regularity conditions, we only 
require that, at the singular point, where we have not a 
single normal direction and a single direction of 
reflection, but a 
cone of normal directions and a cone of directions of 
reflection, there is a direction of reflection that points 
strictly inside the approximating cone, and a normal vector that forms an angle 
of strictly less than $\pi /2$ with every direction of 
reflection (Conditions \ref{G} $(iii)$ and $(iv)$). Conditions 
\ref{G} $(iii)$ and $(iv)$ are the analog of the well known 
completely-{\em S\/} condition for an orthant with constant 
direction of reflection on each face. 

Our argument follows the general outline of 
\cite{KW91}, with two fundamental changes: We 
focus on semimartingale reflecting diffusions and 
characterize them as solutions of 
{\em constrained martingale problems\/} rather than of 
submartingale problems; We replace the Krein-Rutman 
theorem by our reverse ergodic 
theorem for inhomogeneous killed Markov chains. 
Preliminarly we apply 
the Markov selection results of \cite{CK19}, so that we 
can reduce to proving 
uniqueness among strong Markov reflecting diffusions. 
The characterization of a semimartingale reflecting diffusion as a 
solution of a constrained martingale problem (more 
precisely a natural solution) is equivalent to the 
characterization as a solution of a stochastic differential 
equation used in \cite{TW93} and 
\cite{DW96}, but it avoids the need of the 
oscillation estimates required there. 
In order to apply our 
ergodic theorem, we need to prove our uniform lower 
bounds: $(i)$ and $(ii)$ in Theorem \ref{th:ergodic}. We obtain 
the bound $(i)$ by means of some auxiliary functions that 
we construct by elaborating on the functions $\psi_{\alpha}$ and 
$\chi$ introduced in \cite{KW91} (or the corresponding 
functions introduced in \cite{VW85}, in the 2-dimensional 
case). In order to do this, we have to prove that, under 
our conditions, the parameter $\alpha$ of \cite{KW91} (the one whose sign 
determines whether the reflecting Brownian motion 
hits the vertex of the cone) is strictly less than $1$ (Theorem \ref{th:alpha}). 
In fact, under our conditions, we can directly 
construct the reflecting Brownian 
motion as a semimartingale (Theorem \ref{th:cone-exist};) By 
modifying the proof of Theorem 5 of \cite{Wil85}, we show 
that this implies that $\alpha <1$.  
Thus we partially extend the 2-dimensional results of 
\cite{Wil85} to higher dimension. 
In order to verify the bound $(ii)$, we use a coupling 
argument based on Lemma 5.3 of \cite{CK18} (Lemma 
\ref{th:coupling}), the fact 
that, for any reflecting diffusion, $X$, 
the rescaled process $2^{2n}\,X(2^{-4n}\cdot )$ 
converges to a reflecting Brownian motion in the cone 
(the already mentioned Lemma \ref{th:scaling}), and the 
support theorem of \cite{KW91}. 

A recent paper (\cite{Cos23}) shows 
that, under suitable assumptions,  
the results of Section \ref{sectiononesing} of this paper 
extend to obliquely reflecting 
Brownian motion in a nonpolyhedral, piecewise smooth 
cone. The uniqueness result of \cite{Cos23} follows by applying 
Theorem \ref{th:ergodic}. The proofs of the lemmas in 
Subsections \ref{sectionhittimes} and \ref{sectionhitdist}, 
needed to verify the assumptions of Theorem 
\ref{th:ergodic}, carry over essentially without any 
modification, except for Lemma \ref{th:scaling}, which 
requires an additional argument. 
The existence result of \cite{Cos23} follows by modifying 
the proof of Theorem \ref{th:cone-exist}. 
As it is well known, in at least two 
notable classes of stochastic networks, 
namely bandwidth sharing models under a weighted $\alpha$-fair 
bandwidth sharing policy (\cite{KW04} and \cite{KKLW09}) 
and input-queued switches 
under a maximum weight-$\alpha$ policy (\cite{SW12} and \cite{KW12}), 
for all values different from $1$ of the parameter 
$\alpha$, the diffusion approximation for the workload is 
conjectured to be an obliquely reflecting Brownian 
motion in a nonpolyhedral, piecewise smooth cone 
(note that $\alpha$ here is not the parameter of \cite{KW91}.) 
The conjecture for the $\alpha\neq 1$ case 
could not be proved in particular due to the lack of a 
unique characterization for the conjectured limiting 
process. Now, at least in the $\alpha\geq 2$ case,  
\cite{Cos23} provides conditions under which 
the stochastic differential equation that should 
characterize the limiting process has a unique solution. 
In \cite{Cos23} these conditions are actually verified in a 
specific example of bandwidth sharing model. 
Some preliminary computations indicate that a refinement to the $\alpha >1$ case should be possible. 

In dimension $2$, any piecewise ${\cal C}^1$ domain, near each 
singular point that is not a cusp, looks like the domain studied in Section 
\ref{sectiononesing}. In fact the results of Section \ref{sectiononesing}, 
combined with a localization result for constrained 
martingale problems, and with \cite{CK18}, yield new and optimal conditions for existence and 
uniqueness of obliquely reflecting diffusions 
in piecewise ${\cal C}^1$ domains in dimension $2$ (\cite{CK24}). 
These new conditions are strictly more general and more easily 
verifiable than the previously known conditions of 
\cite{DI93}.   
Obliquely reflecting diffusions in $2$-dimensional piecewise smooth domains 
arise naturally, for instance, in singular stochastic control problems. 
As an example, in the problem studied in \cite{WCM94}, 
an optimal control is found from a 
reflecting Brownian motion in a certain piecewise smooth 
domain, if such a process exists. In this example, when the dimension is $2$, the conditions of 
\cite{CK24} seem to be satisfied, but some technical details 
still need to be verified. (See \cite{DF23} for a discussion and a recent contribution on the issue of characterizing optimal singular stochastic controls.)

The outline of the paper is the following: In Section 
\ref{sectionergodic} we prove our reverse 
ergodic theorem for inhomogeneous killed Markov chains, 
while in Section \ref{sectiononesing} we 
prove our existence and uniqueness result for semimartingale reflecting diffusions. 
Section \ref{sectiononesing} is divided into several 
subsections: In Subsection \ref{sectionassum}, we 
state the assumptions, formulate the constrained martingale problem 
and prove some preliminary results; In Subsection 
\ref{sectionexist}, we prove 
existence of a strong Markov, semimartingale reflecting 
diffusion; In Subsection \ref{sectionuniq}, we show how 
uniqueness follows from the reverse ergodic theorem for 
inhomogeneous killed Markov chains; 
In Subsections \ref{sectionhittimes} and 
\ref{sectionhitdist}, we prove the bounds required to apply the reverse ergodic theorem. Finally, 
in Appendix \ref{appcone}, we summarize the results of 
\cite{KW91}, \cite{VW85} and \cite{Wil85} and we prove 
our new results for the cone; Appendix \ref{appaux} contains 
the constructions of various auxiliary functions. 

\subsection{Notation}\label{sectionnotation}
For a set $E$ (in a topological space), we denote by $\stackrel {
\circ}E$ the 
interior of $E$ and by $\overline E$ its closure. 

For a metric space $E$, we denote by ${\cal B}(E)$ the $\sigma$-algebra of 
Borel subsets of $E$ and by ${\cal P}(E)$ and ${\cal M}(E)$ the sets of 
probability measures and of finite signed measures on 
$(E,{\cal B}(E))$ respectively. On ${\cal M}(E)$ we define the total 
variation norm by $||\mu ||_{TV}:=\sup_{C\in {\cal B}(E)}|\mu (C)
|$. For $z\in E$, $\delta_z$ 
denotes the Dirac measure. 
We denote by ${\cal C}(E)$ the space of continuous functions on $
E$, endowed 
with the $\sup$ norm. $C_E[0,\infty )$ is the space of continuous 
functions from $[0,\infty )$ to $E$ and $D_E[0,\infty )$ is the space of 
right continuous functions with left hand limits from $[0,\infty 
)$ to $E$. 

For a stochastic process $Z$, $\{{\cal F}^Z_t\}$ is the 
filtration generated by $Z$, that is ${\cal F}^Z_t:=\sigma (Z(s),\,
s\leq t)$, $t\geq 0$, 
and we use the notation 
\[\P^z(\cdot ):=\P (\cdot |Z(0)=z),\quad\E^z[\cdot ]:=\E[\cdot |Z
(0)=z].\]

In $\R^d$, we denote by $B_r(0)$  the ball of radius $r$ centered at the origin 
and by $S^{d-1}$ the unit sphere. 

\setcounter{equation}{0}

\section{A reverse ergodic 
theorem for inhomogeneous killed Markov chains}\label{sectionergodic}

Let $E$ be a compact metric space, and let $Q(x,dy)$ be a sub-probability transition 
function on $E$, that is, for each $x\in E$, $Q(x,\cdot)$ is a finite measure on $E$ 
with $Q(x,E)\leq 1$ and for each $C\in {\cal B}(E)$, $Q(x,C)$ is 
Borel measurable in $x$.  We will still denote by $Q$ the integral 
operator defined by $Q$. In the proof of their Theorem 5.4,  
\cite{BP87} show that, if $Q$ satisfies the conditions of 
the Krein-Rutman theorem (Theorems 6.1 and 6.3 of 
\cite{KR50}), then, for all $f,g\in C(E)$, $g>0$ 
\footnote{The conditions of the Krein-Rutman theorem actually 
allow for $g\geq 0$ as long as $g$ is not identically zero, but $
g>0$ 
is enough for the application in \cite{BP87} and for our 
purposes as well.}, $\{\mu_k\}\subset {\cal P}(E)$, 
\[\lim_{k\rightarrow\infty}\frac {\int Q^kf(x)\mu_k(dx)}{\int Q^k
g(x)\mu_k(dx)}=C(f,g),\]
where the constant $C(f,g)$ is independent of $\{\mu_k\}$.

Our goal is to extend this result to a sequence of 
compact metric spaces $E_0,E_1,E_2,\ldots$ and a sequence of 
sub-probability transition functions $Q_1,Q_2,\ldots$, with $Q_l$ governing 
transitions from $E_l$ to $E_{l-1}$, 
and give conditions under which $(\ref{eq:constlim})$, i.e. 
\[\lim_{k\rightarrow\infty}\frac {\int Q_kQ_{k-1}\cdots Q_1f(x)\mu_
k(dx)}{\int Q_kQ_{k-1}\cdots Q_1g(x)\mu_k(dx)}=C(f,g)\]
with $C(f,g)$ independent of $\{\mu_k\}$, holds. We may as well take 
$g\equiv 1$, and we will do so in the sequel.

\begin{lemma}\label{th:Ts}
Assume  
\begin{equation}\inf_{x\in E_l}Q_l(x,E_{l-1})>0,\quad\forall l,\label{eq:mass}\end{equation}
and set, for $f\in {\cal C}(E_0)$,  
\begin{equation}T_kf(x):=\frac {Q_kQ_{k-1}\cdots Q_1f(x)}{Q_kQ_{k
-1}\cdots Q_11(x)}.\label{eq:Ts}\end{equation}
If there exists a constant $C(f)$ such that 
\[\sup_{x\in E_k}|T_kf(x)-C(f)|\rightarrow_{k\rightarrow\infty}0,\]
then (\ref{eq:constlim}) holds for $f$ and $g=1$. 
\end{lemma}

\begin{proof}
Divide and multiply by $Q_kQ_{k-1}\cdots Q_11(x)$ inside the integral 
in the numerator of (\ref{eq:constlim}) (with $g=1$).
\end{proof}

\begin{remark}\label{re:subprobseq}
Lemma \ref{th:Ts}, and hence Theorem \ref{th:ergodic} 
below, hold for a sequence of nonzero subprobability measures 
$\{\mu_k\}$ as well. 
\end{remark}
\vskip.2in

Note that the operator $T_k$ defined by $(\ref{eq:Ts})$ 
corresponds to a probability transition function from $E_k$ to $E_
0$
and can be written as 
\begin{equation}T_kf(x)=P_kP_{k-1}\cdots P_1f(x)\label{eq:norm-op}\end{equation}
where the $P_l$ are the operators corresponding to 
the probability transition functions from $E_l$ into $E_{l-1}$ given by
\begin{eqnarray}\label{eq:norm-ker}
P_1(x,dy):=\frac {Q_1(x,dy)}{Q_11(x)}=\frac {Q_1(x,dy)}{\int_{E_
0}Q_1(x,dz)}\qquad\qquad\qquad\qquad\non\\
\quad\\
P_l(x,dy):=\frac {Q_l(x,dy)[Q_{l-1}\cdots Q_11(y)]}{Q_l\cdots Q_
11(x)}=\frac {Q_l(x,dy)[Q_{l-1}\cdots Q_11(y)]}{\int_{E_{l-1}}Q_{
l-1}\cdots Q_11(z)Q_l(x,dz)}\quad l\geq 2.\non\end{eqnarray}

\begin{lemma}\label{th:tvbnd}
Assume $(\ref{eq:mass})$. Define 
\begin{equation}f_{l,\tilde {x}}(x,y):=\frac {dQ_l(x,\cdot )}{d\big(Q_l(x,\cdot 
)+Q_l(\tilde {x},\cdot )\big)}(y)\label{eq:reldens}\end{equation}
and
\begin{equation}\epsilon_l(x,\tilde {x}):=\int\left(f_{l,\tilde {x}}(x,y)\wedge 
f_{l,x}(\tilde x,y)\right)(Q_l(x,dy)+Q_l(\tilde {x},dy)).\label{eq:tvdensbnd}\end{equation}
Then, for $P_l$ given by $(\ref{eq:norm-ker})$, 
\begin{eqnarray}
\Vert P_l(x,\cdot )-P_l(\tilde {x},\cdot )\Vert_{TV}&\leq&1-\epsilon_
l(x,\tilde {x})\inf_{z\in E_l,y\in E_{l-1}}\left(\frac {Q_{l-1}\cdots 
Q_11(y)}{Q_l\ldots Q_11(z)}\right)\label{eq:tvbnd}\\
&\leq&1-\epsilon_l(x,\tilde {x})\inf_{z,y\in E_{l-1}}\left
(\frac {Q_{l-1}\cdots Q_11(y)}{Q_{l-1}\ldots Q_11(z)}\right)\nonumber\end{eqnarray}
\end{lemma}

\begin{proof}
Observe that $P_l(x,dy)<<Q_l(x,dy)$ with density given by 
$(\ref{eq:norm-ker})$. Then  
\begin{eqnarray*}
&&\Vert P_l(x,\cdot )-P_l(\tilde {x},\cdot )\Vert_{TV}\\
&&=\frac 12\int\left|f_{l,\tilde {x}}(x,y)\frac {Q_{l-1}\cdots 
Q_11(y)}{Q_l\ldots Q_11(x)}-f_{l,x}(\tilde x,y)\frac {Q_{l-1}\cdots 
Q_11(y)}{Q_l\ldots Q_11(\tilde {x})}\right|(Q_l(x,dy)+Q_l(\tilde {
x},dy))\\
&&=1-\int\hskip-.03in\bigg(f_{l,\tilde {x}}(x,y)\frac {Q_{l-1}\cdots Q_
11(y)}{Q_l\ldots Q_11(x)}\wedge f_{l,x}(\tilde x,y)\frac {Q_{l-1}
\cdots Q_11(y)}{Q_l\ldots Q_11(\tilde {x})}\bigg)(Q_l(x,dy)+Q_l(
\tilde {x},dy))\\
&&\leq 1-\int\hskip-.03in\big(f_{l,\tilde {x}}(x,y)\wedge f_{l,x}(\tilde 
x,y)\big)\bigg(\frac {Q_{l-1}\cdots Q_11(y)}{Q_l\ldots Q_11(x)\vee 
Q_l\ldots Q_11(\tilde {x})}\bigg)(Q_l(x,dy)+Q_l(\tilde {x},dy))\\
&&\leq 1-\epsilon_l(x,\tilde {x})\inf_{z\in E_l,y\in E_{l-1}}\left
(\frac {Q_{l-1}\cdots Q_11(y)}{Q_l\ldots Q_11(z)}\right)\end{eqnarray*}
The second inequality in $(\ref{eq:tvbnd})$ follows from the fact that
\[Q_l\ldots Q_11(x)\leq Q_l(x,E_{l-1})\sup_{z\in E_{l-1}}Q_{l-1}\cdots 
Q_11(z)\leq\sup_{z\in E_{l-1}}Q_{l-1}\cdots Q_11(z).\]
\end{proof}

\begin{theorem}\label{th:ergodic}
For $x,\tilde {x}\in E_l$, let $\epsilon_l(x,\tilde {x})$ be defined as in Lemma 
\ref{th:tvbnd}. Assume $Q_l$ is not identically zero, i.e. 
$\sup_xQ_l(x,E_{l-1})>0$, for all $l$, and there exist $c_0>0$ and $
\epsilon_0>0$ such 
that 

\begin{itemize}
\item[(i)] 
\[\inf_{x\in E_k}Q_k\cdots Q_11(x)\geq c_0\sup_{x\in E_k}Q_k\cdots 
Q_11(x),\quad\forall k,\]
\item[(ii)] 
\[\inf_k\inf_{x,\tilde {x}\in E_k}\epsilon_k(x,\tilde {x})\geq\epsilon_
0.\]
\end{itemize}
Then 
\begin{equation}\inf_{x\in E_k}Q_k\cdots Q_11(x)>0,\quad\forall k
,\label{Qn0}\end{equation}
and $(\ref{eq:constlim})$ holds for all $f\in {\cal C}(E_0)$ and $
g=1$.  
\end{theorem}

\begin{proof}
First of all note that (i) above and the assumption that, 
for every $l$, $Q_l$ is never identically zero imply, by 
induction, $(\ref{Qn0})$, 
which, in turn, implies $(\ref{eq:mass})$. 

Next, for $\mu\in {\cal P}(E_l)$, denote 
\[\mu P_l(dy):=\int_{E_l}P_l(x,dy)\mu (dx).\]
Of course we can suppose $\epsilon_0<1$, $c_0<1$. 
For $\mu ,\tilde{\mu}\in {\cal P}(E_l)$, by Lemma \ref{th:tvbnd}, for all $
l$,
\begin{eqnarray*}
&&\|\mu P_l-\tilde{\mu }P_l\|_{TV}=\sup_{C\in {\cal B}(E_{l-1})}\bigg
|\int_{E_l}\int_{E_l}\big(P_l(x,C)-P_l(\tilde {x},C))\big)\mu (dx
)\tilde{\mu }(d\tilde {x})\bigg|\\
&&\leq\sup_{x,\tilde {x}\in E_l}\|P_l(x,\cdot )-P_l(\tilde {x},\cdot 
)\|_{TV}\leq 1-\epsilon_0c_0.\end{eqnarray*}
Then, by Lemma 5.4 of \cite{CK18}, for $\mu ,\tilde{\mu}\in {\cal P}
(E_k)$,
\begin{eqnarray*}
&&\|\mu P_kP_{k-1}\cdots P_1-\,\tilde{\mu }P_kP_{k-1}\cdots P_1\|_{
TV}=\|(\mu P_k)(P_{k-1}\cdots P_1)-\,(\tilde{\mu }P_k)(P_{k-1}\cdots 
P_1)\|_{TV}\\
&&\leq\|\mu P_k-\tilde{\mu }P_k\|_{TV}\|\mu P_{k-1}\cdots P_1-\tilde{
\mu }P_{k-1}\cdots P_1\|_{TV}\\
&&\leq (1-\epsilon_0c_0)\|\mu P_{k-1}\cdots P_1-\,\tilde{\mu }P_{
k-1}\cdots P_1\|_{TV},\end{eqnarray*}
and, by iterating,
\[\|\mu P_kP_{k-1}\cdots P_1-\,\tilde{\mu }P_kP_{k-1}\cdots P_1\|_{
TV}\leq (1-\epsilon_0c_0)^k.\]
In particular, for each $f\in {\cal C}(E_0)$, for an arbitrary $\{
x_k\}$, $x_k\in E_k$ for 
each $k$,  
\begin{eqnarray*}|T_{k+l}f(x_{k+l})-T_kf(x_k)|&\leq&\|(\delta_{x_{k+l}}P_{k+l}\cdots 
P_{k+1})P_k\cdots P_1-\delta_{x_k}P_k\cdots P_1\|_{TV}\|f\|\\
&\leq& (1-\epsilon_0c_0)^k\|f\|,\end{eqnarray*}
so that $\{T_kf(x_k)\}$ is a Cauchy sequence. If $C(f)$ is its 
limit, we get, in an analogous manner,  
\[\sup_{x\in E_k}|T_kf(x)-C(f)|\leq (1-\epsilon_0c_0)^k\|f\|+|T_k
f(x_k)-C(f)|,\]
which yields the assertion by Lemma \ref{th:Ts}. 
\end{proof}

\section{Existence and uniqueness for reflecting diffusions in a domain with one singular point}\label{sectiononesing}

\setcounter{equation}{0}

\subsection{Formulation of the problem and preliminaries} 
\label{sectionassum}

We consider a connected domain $D\subseteq\R^d$ that has a smooth 
boundary except at a 
single point, which we will take to be the origin, and 
that in a neighborhood of the singular point can be 
approximated by a 
cone. More precisely we assume the following condition 
$(d_H$ denotes the Hausdorff distance). 

\begin{condition}\label{D}\hfill 
\begin{itemize}
\item[(i)]$D$ is a bounded domain and $\partial D-\{0\}$ is of class 
${\cal C}^1$.  There exist a nonempty 
domain, ${\cal S}$, in the unit sphere, $S^{d-1}$, $r_D>0$ and $c_
D>0$ such that, setting 
\[{\cal K}:=\{x:\,x=rz,\,z\in {\cal S},\,r>0\},\]
for $r\leq r_D$, 
\begin{equation}d_H(\overline D\cap\partial B_r(0),\overline {{\cal K}}
\cap\partial B_r(0))\leq c_Dr^2,\label{eq:condi1}\end{equation}
\begin{equation}\sup_{x\in\partial D\cap\partial B_r(0)}d(x,\partial 
{\cal K}\cap\partial B_r(0))\leq c_Dr^2.\label{eq:condi2}\end{equation}
\item[(ii)]For $x\in\partial D-\{0\}$ denote by $n(x)$ the unit inward 
normal to $\overline D$ and for $x\in\partial {\cal K}-\{0\}$ denote by $
n^{{\cal K}}(x)$ the unit inward 
normal to $\overline {{\cal K}}$. Note that $n^{{\cal K}}$ has the scaling property: $
n^{{\cal K}}(rz)=n^{{\cal K}}(z)$ 
for all $z\in\partial {\cal S}$, $r>0$. 
For $x\in\partial D-\{0\}$, $|x|\leq r_D$, and $\bar {z}\in\partial 
{\cal S}$ such that 
$|\frac x{|x|}-\bar {z}|=d(\frac x{|x|},\partial {\cal S})$, 
\begin{equation}|n(x)-n^{{\cal K}}(\bar {z})|\leq c_D|x|.\label{eq:condii}\end{equation}
\item[(iii)]For $d\geq 3$, the boundary $\partial {\cal S}$ of ${\cal S}$ in $
S^{d-1}$ is of class ${\cal C}^3$.
\end{itemize}
\end{condition}

\begin{remark}\label{re:D}
Clearly $(\ref{eq:condi2})$ implies that, for 
every $x\in\partial D$, $|x|\leq r_D$, 
\[d(\frac x{|x|},\partial {\cal S})\leq c_{\!D}|x|,\]
and hence, by the smoothness of $\partial {\cal S}$ (Condition \ref{D}(iii)), 
for $r_D$ sufficiently small and $|x|\leq r_D$, there is a unique $
\bar {z}\in\partial {\cal S}$ such 
that 
\[|\frac x{|x|}-\bar {z}|=d(\frac x{|x|},\partial {\cal S}).\]

The same for $x\in\overline D-\overline {{\cal K}}$.
\end{remark}
\vskip.2in

We assume the following on the directions of reflection.

\begin{condition}\label{G}\hfill
\begin{itemize}
\item[(i)]$g:\R^d-\{0\}\rightarrow\R^d$ is a locally Lipschitz continuous vector 
field, of unit length on $\partial D-\{0\}$, such that 
\[\inf_{x\in\partial D-\{0\}}g(x)\cdot n(x)>0.\]
There exists a unit vector field $\bar {g}:\partial {\cal S}\rightarrow\R^
d$ and $c_g>0$ such 
that, for $x\in\partial D-\{0\}$ and $\bar {z}$ the closest point to $\frac 
x{|x|}$ 
on $\partial {\cal S}$, 
\[|g(x)-\bar {g}(\bar {z})|\leq c_g|x|,\quad\mbox{\rm for }|x|\leq 
r_D.\]

\item[(ii)]Extending the definition of $\bar {g}$ to $\partial {\cal K}
-\{0\}$ by
\[\bar {g}(x):=\bar {g}(\frac x{|x|}),\]
\[\inf_{x\in\partial {\cal K}-\{0\}}\bar {g}(x)\cdot n^{{\cal K}}
(x)>0.\]

For $d\geq 3$, $\bar {g}$ is of class ${\cal C}^2$ and, for $x\in
\partial {\cal K}-\{0\}$, denoting by 
$v^r(x):=\frac x{|x|}$,  the radial unit vector, 
$\frac {\bar {g}\cdot v^r}{\bar {g}\cdot n^{{\cal K}}}$ is of class $
{\cal C}^3$ and $\frac {\bar {g}-\bar {g}\cdot v^rv^r-\bar {g}\cdot 
n^{\!{\cal K}}n^{{\cal K}}}{\bar {g}\cdot n^{{\cal K}}}$ is of class $
{\cal C}^4$ 
(see Remark \ref{re:KW}.)

\item[(iii)]For $x\in\partial D-\{0\}$, let $G(x):=\{\eta g(x),\,
\eta\geq 0\}$, and let 
$G(0)$ be the closed, convex cone generated by 
$\{\bar {g}(z),\,z\in\partial {\cal S}\}$. Assume 
\[G(0)\cap {\cal K}\neq\emptyset .\]

\item[(iv)]Let $N(0)$ denote the normal cone at the origin 
for $\overline D$, that is,
\begin{equation}N(0):=\{n\in\R^d:\,\liminf_{x\in\overline D-\{0\},\,x\rightarrow 
0}n\cdot\frac x{|x|}\geq 0\,\}.\label{eq:Ncone}\end{equation}
There exists a unit vector $e\in N(0)$ such that
\begin{equation}\inf_{g\in G(0),\,|g|=1}e\cdot g=c_e>0,\label{eq:condGiv1}\end{equation}
and, possibly by taking a smaller $r_D$, 
\begin{equation}e\cdot x\geq 0,\quad\forall x\in\overline D,\;|x|
\leq r_D.\label{eq:condGiv2}\end{equation}
$(\ref{eq:condGiv2})$ is always satisfied if $D$ is convex or 
if $\stackrel {\circ}{N(0)}\neq\emptyset$: See Remark \ref{re:Ncone} below.
\end{itemize}
\end{condition}

\begin{remark}\label{re:Ncone}
$N(0)$ is a closed convex cone. 
By Condition \ref{D}(i), $N(0)$ is also the normal cone at 
the origin for $\overline {{\cal K}}$, i.e. 
\begin{equation}N(0)=\{n\in\R^d:\,n\cdot x\geq 0,\,\forall x\in\overline {{\cal K}}
\}.\label{eq:KNcone}\end{equation}
Condition \ref{G}(iv) implies that 
$N(0)\neq\emptyset$, hence $\overline {{\cal K}}$ is always contained in a closed halfspace. 

If $\stackrel {\circ}{N(0)}$ is nonempty, we can assume, 
without loss of generality, that $e\in\stackrel {\circ}{N(0)}$. Then $\overline {
{\cal K}}-\{0\}$ is 
contained in the open half space $\{x\in\R^d:\,x\cdot e>0\}$, and 
hence, by $(\ref{eq:condi1})$, possibly taking a smaller $r_D$, 
$\big(\overline D-\{0\}\big)\cap\overline {B_{r_D}(0)}$ is contained in the same open half 
space. 

Assuming $(\ref{eq:condGiv2})$
instead of $\stackrel {\circ}{N(0)}\neq\emptyset$ allows, for instance, to deal with 
same cases in which $\partial D$ is actually 
smooth but $g$ is discontinuous at the origin.
$(\ref{eq:condGiv2})$ is used only in the proof of 
Lemma \ref{th:extime}. 
\end{remark}
\vskip.2in

\begin{remark}\label{re:KW}
Conditions \ref{D}(iii) and \ref{G}(ii) are the assumptions 
of \cite{KW91}, which we need because we will exploit 
some of their results. 
\end{remark}
\vskip.2in

Reflecting diffusions are often characterized as 
solutions of stochastic differential equations. Assume the 
following. 

\begin{condition}\label{A}\hfill
\begin{itemize}
\item[(i)] $b:\R^d\rightarrow\R^d$ and $\sigma :\R^d\rightarrow\R^{
d\times d}$ are Lipschitz continuous. 
\item[(ii)] $\sigma (0)$ is non singular. 
\end{itemize}
\end{condition}

\begin{definition}\label{def:SDER}
A continuous process $X$ is a solution of the {\em stochastic }
{\em differential equation with reflection\/} in $\overline D$ with coefficients $
b$ and $\sigma$ 
and direction of reflection $g$, if there exist a 
standard Brownian motion $W$, a continuous, non 
decreasing process $\lambda$ and a process 
$\gamma$ with measurable paths, all defined on the same probability space as $
X$, 
such that $W(t+\cdot )-W(t)$ is independent of 
${\cal F}_t^{X,W,\lambda ,\gamma}$ for all $t\geq 0$ and  
\begin{eqnarray}
&X(t)=X(0)+\int_0^tb(X(s))ds+\int_0^t\sigma (X(s))dW(s)+\int_0^t\gamma 
(s)\,d\lambda (s),\quad t\geq 0,&\non\\
&\quad\gamma (t)\in G(X(t)),\quad |\gamma (t)|=1,\quad d\lambda -
a.e.,\quad t\geq 0,&\label{eq:SDER}\\
&X(t)\in\overline D,\quad\lambda (t)=\int_0^t{\bf 1}_{\partial D}
(X(s))d\lambda (s),\quad t\geq 0,&\non\end{eqnarray}
is satisfied a.s.. 

Given an initial distribution $\nu\in {\cal P}(\overline D)$, {\em weak }
{\em uniqueness\/} or {\em uniqueness in distribution\/} holds if any 
two solutions of (\ref{eq:SDER}) with
$P\{X(0)\in\cdot \}=\nu$ have the same distribution on 
$C_{\overline D}[0,\infty )$. 

A stochastic process $\tilde {X}$ (for example a solution of an 
appropriate martingale problem or submartingale problem) 
is a weak solution of $(\ref{eq:SDER})$ if there is a 
solution $X$ of (\ref{eq:SDER}) such that $\tilde {X}$ and $X$ have the 
same distribution.

We say that $(\ref{eq:SDER})$ admits a strong solution 
if, given a filtered probability space $(\Omega ,{\cal F},$ $\{{\cal F}_
t\},\P )$, a 
standard $\{{\cal F}_t\}$-Brownian motion $W$ and a ${\cal F}_0$-measurable, $\overline 
D$-valued random 
variable $X_0$, there exist a continuous  
process $X$, with $X(0)=X_0$, 
a continuous, non decreasing process $\lambda$ and a process 
$\gamma$ with measurable paths, all $\{{\cal F}_t\}$-adapted, such that $
(\ref{eq:SDER})$ is 
satisfied a.s.. 

We say that strong (or pathwise) uniqueness holds for $(\ref{eq:SDER})$ 
if any two solutions defined on the same probability 
space with the same Brownian motion and the same 
initial condition coincide for all times a.s..
\end{definition}

We denote by $A$ the operator 
\begin{equation}{\cal D}(A):={\cal C}^2(\overline D),\quad Af(x):
=b(x)\cdot\nabla f(x)+\frac 12\mbox{\rm tr}((\sigma\sigma^T)(x)D^
2f(x)).\label{eq:A}\end{equation}

\begin{remark}\label{re:sigma(0)}
Since $\sigma (0)$ is non singular, it is easy to check that $X$ 
is a solution to $(\ref{eq:SDER})$ if and only if $\sigma^{-1}(0)
X$ is a solution of 
$(\ref{eq:SDER})$, with the appropriate coefficients, in the 
corresponding domain $\sigma^{-1}(0)D$, with 
direction of reflection $\frac {\sigma^{-1}(0)g\circ\sigma (0)}{|
\sigma^{-1}(0)g\circ\sigma (0)|}$, and the new domain and 
direction of reflection satisfy Condition 
\ref{D} and Condition \ref{G} if and only if $D$ and $g$ do. 
Therefore, without loss of 
generality, we will take, from now on, 
\[\sigma (0)=I.\]
\end{remark}
\vskip.2in

Conditions \ref{D}, \ref{G} and \ref{A} will be our 
standing assumptions. 
\vskip.2in

Reflecting diffusions can also be characterized 
as {\em natural solutions\/} of {\em constrained martingale problems}. 
Constrained martingale problems have 
been introduced by \cite{Kur90} and \cite{Kur91}) for 
operators on functions defined on metric spaces, and 
further studied in \cite{KS01}, \cite{CK15} and \cite{CK19}.  
Here we consider the specific constrained 
martingale problem corresponding to $(\ref{eq:SDER})$, 
which is formulated as follows. 

Let $G(x)$, $x\in\partial D-\{0\}$, and $G(0)$ be as in Condition \ref{G} (iii), and let 
\begin{eqnarray}
&&\qquad\qquad\qquad\qquad U:=S^{d-1},\non\\
&&\qquad\Xi :=\{(x,u)\in\partial D\times U:u\in G(x)\},\label{eq:clset}\\
&&B:{\cal D}(B):={\cal C}^2(\overline D)\,\rightarrow\,{\cal C}(\Xi 
),\quad Bf(x,u):=\nabla f(x)\cdot u\non\end{eqnarray}
The following argument shows that $\Xi$ is closed: Suppose 
$\{(x^n,u^n)\}\subseteq\Xi$, $(x^n,u^n)\rightarrow (x,u)$ as $n\rightarrow
\infty$; If $x\neq 0$, eventually 
$x^n\neq 0$, $u^n=g(x^n)$ and $u=g(x)$ by the continuity of $g$; If 
$x=0$ and there exists a subsequence $\{x^{n_i}\}$ such that 
$x^{n_i}=0$ for all $i$, then $u^{n_i}\in G(0)$ for all $i$, hence $
u\in G(0)$ 
because $G(0)$ is closed; If eventually $x^n\neq 0$, then 
eventually $|g(x^n)-\bar {g}(\bar {z}^n)|\leq c_g|x^n|$, where $\bar {
z}^n$ is the closest 
point to $x^n/|x^n|$ on $\partial {\cal S}$, so that $\bar {g}(\bar {
z}^n)\rightarrow u$ as $n\rightarrow\infty$ and 
again $u\in G(0)$ because $G(0)$ is closed.

Define ${\cal L}_U$ to be the space of measures $\mu$ on $[0,\infty 
)\times U$ 
such that $\mu ([0,t]\times U)<\infty$ for all $t>0$.  ${\cal L}_
U$ 
is topologized so that $\mu_n\in {\cal L}_U\rightarrow\mu\in {\cal L}_
U$ if and only if 
\[\int_{[0,\infty )\times U}f(s,u)\mu_n(ds\times du)\rightarrow\int_{
[0,\infty )\times U}f(s,u)\mu (ds\times du)\]
for all continuous $f$ with compact support in $[0,\infty )\times 
U$. It is 
possible to define a metric on ${\cal L}_U$ that induces the above 
topology and makes ${\cal L}_U$ into a complete, separable metric 
space. Also define ${\cal L}_{\Xi}$ in the same way. 

We will say that an ${\cal L}_U$-valued random variable 
$\Lambda_1$ is adapted to a filtration $\{{\cal G}_t\}$ if 
\[\Lambda_1([0,\cdot ]\times C)\mbox{\rm \ is }\{{\cal G}_t\}-\mbox{\rm adapted}
,\quad\forall C\in {\cal B}(U).\]
We define an adapted ${\cal L}_{\Xi}$-valued random variable 
analogously. 

\begin{definition}\label{def:cmp}
Let $A$, $\Xi$ and $B$ be as in $(\ref{eq:A})$ and $(\ref{eq:clset})$. 
A process $X$ in $D_{\overline D}[0,\infty )$ is a solution of the 
{\em constrained martingale problem\/} for $(A,D,B,\Xi )$ if 
there exists a random measure $\Lambda$ with values in ${\cal L}_{
\Xi}$ and a filtration 
$\{{\cal F}_t\}$ such that $X$ and $\Lambda$ are $\{{\cal F}_t\}$-adapted and for each 
$f\in {\cal C}^2(\overline D)$, 
\begin{equation}f(X(t))-f(X(0))-\int_0^tAf(X(s))ds-\int_{[0,t]\times
\Xi}Bf(x,u)\Lambda (ds\times dx\times du)\label{eq:cmp}\end{equation}
is a $\{{\cal F}_t\}$-local martingale. 
By the continuity of $f$, we may  
assume, without loss of generality, that $\{{\cal F}_t\}$ is right 
continuous. 

Given $\nu\in {\cal P}(\overline D)$, we say that there is a unique solution of 
the constrained martingale problem for $(A,D,B,\Xi )$ with 
initial distribution $\nu$ if any two solutions 
with $P\{X(0)\in\cdot \}=\nu$ have the same distribution on $D_{\overline 
D}[0,\infty )$. 

\end{definition}

\begin{remark}\label{re:semimg}
Since $f(x):=x_i$ $i=1,...,d$, belongs to ${\cal D}(A)={\cal D}(B
)$, every 
solution of the constrained martingale problem 
for $(A,D,B,\Xi )$ is a semimartingale. 
\end{remark}

In general, an effective way of constructing solutions of a 
constrained martingale problem is by time-changing 
solutions of the corresponding {\em controlled martingale }
{\em problem\/} which is a ''slowed down'' version of the 
constrained martingale problem. 

\begin{definition}\label{def:clmp}
Let $A$, $U$, $\Xi$ and $B$ be as in $(\ref{eq:A})$ and $(\ref{eq:clset})$. 
$(Y,\lambda_0,\Lambda_1)$ is a solution of the {\em controlled martingale problem \/}
for $(A,D,B,\Xi )$, if $Y$ is a process in ${\cal D}_{\overline 
D}[0,\infty )$, 
$\lambda_0$ is nonnegative and nondecreasing, 
$\Lambda_1$ is a random measure with values in ${\cal L}_U$ such that 
\begin{equation}\lambda_1(t):=\Lambda_1([0,t]\times U)=\int_{[0,t
]\times U}{\bf 1}_{\Xi}(Y(s),u)\Lambda_1(ds\times du),\label{Lam1supp}\end{equation}
\[\lambda_0(t)+\lambda_1(t)=t,\]
and there exists a filtration $\{{\cal G}_t\}$ such that $Y$, $\lambda_
0$, and $\Lambda_1$ are 
$\{{\cal G}_t\}$-adapted and 
\begin{equation}f(Y(t))-f(Y(0))-\int_0^tAf(Y(s))d\lambda_0(s)-\int_{
[0,t]\times U}Bf(Y(s),u)\Lambda_1(ds\times du)\label{eq:clmgp}\end{equation}
is a $\{{\cal G}_t\}$-martingale for all $f\in {\cal C}^2(\overline 
D)$. 
By the continuity of $f$, we may 
assume, without loss of generality, that $\{{\cal G}_t\}$ is right 
continuous. 
\end{definition}

\begin{remark}\label{re:clmpcont}
It can be easily verified (e.g. by Proposition 3.10.3 of 
\cite{EK86}) that, for every 
solution of the controlled martingale problem 
for $(A,D,B,\Xi )$, $Y$ is continuous.
\end{remark}

\begin{definition}\label{def:nat}
Let $A$, $U$, $\Xi$ and $B$ be as in $(\ref{eq:A})$ and $(\ref{eq:clset})$. 
A solution, $X$, of the constrained martingale problem 
for $(A,D,B,\Xi )$ is called {\em natural\/} if, for some 
solution $(Y,\lambda_0,\Lambda_1)$ of the controlled martingale problem 
for $(A,D,B,\Xi )$ with filtration $\{{\cal G}_t\}$, 
\begin{eqnarray}
\label{eq:tchange}&&\qquad\qquad\qquad X(t)=Y(\lambda_0^{-1}(t)),
\quad {\cal F}_t={\cal G}_{\lambda_0^{-1}(t)},\non\\
&&\Lambda ([0,t]\times C)=\int_{[0,\lambda_0^{-1}(t)]\times U}{\bf 1}_
C(Y(s),u)\Lambda_1(ds\times du),\quad C\in {\cal B}(\Xi ),\\
&&\qquad\mbox{\rm where }\lambda_0^{-1}(t)=\inf\{s:\lambda_0(s)>t
\},\quad t\geq 0.\non\end{eqnarray}
Given $\nu\in {\cal P}(\overline D)$, we say that there is a unique natural solution of 
the constrained martingale problem for $(A,D,B,\Xi )$ with 
initial distribution $\nu$ if any two natural solutions 
with $P\{X(0)\in\cdot \}=\nu$ have the same distribution on $D_{\overline 
D}[0,\infty )$. 
\end{definition}

Given a solution $(Y,\lambda_0,\Lambda_1)$ of the controlled martingale problem 
for $(A,D,B,\Xi )$, the time changed process $X$ defined by   
$(\ref{eq:tchange})$ will not always be 
a solution of the corresponding constrained martingale problem. In fact 
it may be impossible to stop $(\ref{eq:clmgp})$, after 
the time change by $\lambda_0^{-1}$, in such 
a way that the stopped process is a local martingale. 
Conditions under which it is possible are given in 
\cite{CK19}, Corollary 3.9, and
the following lemma guarantees that they are 
satisfied under our standing assumptions. 

\begin{lemma}\label{th:bdryfunc}
There exists a function $F\in {\cal C}^2(\overline D)$ such that 
\[\inf_{\,x\in\partial D}\inf_{g\in G(x),\,|g(x)|=1}\nabla F(x)\cdot 
g:=c_F>0.\]
\end{lemma}

\begin{proof}
See Appendix \ref{appaux}. 
\end{proof}

\begin{proposition}\label{th:tchange}
Let $A$, $U$, $\Xi$ and $B$ be as in $(\ref{eq:A})$ and 
$(\ref{eq:clset})$ and assume Conditions \ref{D}, \ref{G} 
and \ref{A}.

For every solution of the controlled martingale problem 
for $(A,D,B,\Xi )$, $\lambda_0(t)\rightarrow\infty$ as $t\rightarrow
\infty$ almost surely and the time changed process $X$ defined by   
$(\ref{eq:tchange})$ is a natural solution of the 
corresponding constrained martingale problem. 
$(\ref{eq:cmp})$ is a martingale. 
\end{proposition}

\begin{proof}
By Lemma \ref{th:bdryfunc}, Lemma 3.1 of \cite{CK19} and Corollary 3.9 
a) of \cite{CK19}, $\lim_{t\rightarrow\infty}\lambda_0(t)=\infty$ a.s. and, after the time change by $
\lambda_0^{-1}$, 
$(\ref{eq:clmgp})$ is a martingale.  
\end{proof}
\vskip.2in

Theorem \ref{th:equiv} below shows that, under our standing assumptions, 
the two characterizations of a reflecting diffusion in 
$\overline D$ with coefficients $b$ and $\sigma$ and direction of reflection $
g$ as a 
solution of $(\ref{eq:SDER})$ and as a natural solution of 
the constrained martingale problem defined 
by $(\mbox{\rm \ref{eq:A}})$$ $ and $(\ref{eq:clset})$ are equivalent. This result parallels the equivalence of the two characterizations of a reflecting diffusion as as a solution of $(\ref{eq:SDER})$ and as a solution of the corresponding submartingale problem proved in Theorem 1 of \cite{KR17}. 

The proof of Theorem \ref{th:equiv} relies on the 
following lemma. 

\begin{lemma}\label{th:takeoff}
For every solution $(Y,\lambda_0,\Lambda_1)$ of the controlled martingale 
problem for $(A,D,B,\Xi )$, $\lambda_0(t)>0$ for all $t>0$, a.s.. 

Moreover, $\lambda_0$ is strictly increasing, a.s..
\end{lemma}

\begin{proof}
The first assertion follows essentially from Condition 
\ref{G}(iv) and Remark \ref{re:clmpcont}. The proof is analogous to that of Lemma 6.8 of 
\cite{CK19} and Lemma 3.1 of \cite{DW96}. The second 
assertion follows from Lemma 3.4 of \cite{CK19}. 
\end{proof}

\begin{remark}\label{re:cmpcont}
It follows from Remark \ref{re:clmpcont} and Lemma 
\ref{th:takeoff} that every natural solution of the 
constrained martingale problem for $(A,D,B,\Xi )$ is a.s. 
continuous and so is the corresponding process $\Lambda ([0,\cdot 
]\times\Xi )$.
\end{remark}
\vskip.2in

\begin{theorem}\label{th:equiv}
Let $A$, $U$, $\Xi$ and $B$ be as in $(\ref{eq:A})$ and 
$(\ref{eq:clset})$ and assume Conditions \ref{D}, \ref{G} 
and \ref{A}.

Every solution of $(\ref{eq:SDER})$ is a 
natural solution of the constrained martingale problem 
for $(A,D,B,\Xi )$. 

Every natural solution of 
the constrained martingale problem 
for $(A,D,B,\Xi )$ is a weak solution of 
$(\ref{eq:SDER})$.
\end{theorem}

\begin{proof}
The proof is the same as for Theorem 6.12 of 
\cite{CK19}. It relies essentially on Lemma \ref{th:takeoff} 
and Remark \ref{re:cmpcont}. 
\end{proof}

\begin{remark}\label{re:alt}
The first assertion of Proposition \ref{th:tchange} could have been proved, 
alternatively, by 
Lemma \ref{th:takeoff} and Lemmas 3.3 and 3.4 of 
\cite{CK19}. 
\end{remark}
\vskip.2in

We conclude this section with two important properties 
of a natural solution $X$ of the 
constrained martingale problem for $(A,D,B,\Xi )$. For $\delta >0$, define
\begin{equation}\tau^{X,\delta}:=\inf\{\{t\geq 0:\,|X(t)|=\delta 
\}.\label{eq:taudelta}\end{equation}
Whenever there is no risk of confusion, we will omit 
the superscript $X$. 

\begin{lemma}\label{th:extime} 
There exists $\bar{\delta }>0$, $\bar{\delta}\leq r_D$, $\bar {c}
>0$, depending only on the data of 
the problem, such that, for $\delta\leq\bar{\delta}$, for every 
natural solution $X$ of the 
constrained martingale problem for $(A,D,B,\Xi )$ starting at 
$0$, 
\[\E[\tau^{X,\delta}]\leq\bar {c}\delta^2.\]
\end{lemma}

\begin{proof}
The assertion follows essentially from Condition 
\ref{G}(iv). 
The proof is analogous to that of Lemma 4.2 of 
\cite{CK18} and Lemma 6.4 of \cite{TW93}.
\end{proof}

\begin{lemma}\label{th:occtime}
For every natural solution $X$ of the constrained martingale 
problem for $(A,D,B,\Xi )$, 
\[\int_0^{\infty}{\bf 1}_{\{0\}}(X(t))\,dt=0,\quad a.s..\]
\end{lemma}

\begin{proof}
The proof uses the same argument as Lemma 2.1 of \cite{TW93}. 
Fix an arbitrary unit vector $v$. Then, by Remark 
\ref{re:cmpcont}, 
\[m(t):=v\cdot X(t)-v\cdot X(0)-\int_0^tv\cdot b(X(s))\,ds\,-\int_{
[0,t]\times\Xi}v\cdot u\,\Lambda (ds\times du)\]
is a continuous semimartingale with 
\[[m,m](t)=\int_0^t|\sigma (X(s))^Tv|^2ds.\]
Therefore, by Tanaka's formula (see, e.g., \cite{Prot90}, Corollary 
1 to Theorem 51, Chapter IV, Section 5), for each $t>0$, 
\[\int_0^t{\bf 1}_{\{0\}}(X(s)\cdot v)\,|\sigma (X(s))^Tv|^2ds=\int_{\R}
{\bf 1}_{\{0\}}(a)\,L_t(a)da=0,\quad a.s.,\]
$L_t(a)$ being the local time of $m$ at $a$. Hence the set of 
times 
$\{s\leq t:\,X(s)\cdot v=0\mbox{\rm \ and }|\sigma (X(s))^Tv|\neq 
0\}$ has zero Lebesgue 
measure, a.s., which yields the assertion by Condition 
\ref{A} (ii). 
\end{proof}

\subsection{Existence} \label{sectionexist}

In this subsection we show that there exists a strong 
Markov, natural solution of 
the constrained martingale problem for $(A,D,B,\Xi )$, and 
hence, by Theorem \ref{th:equiv}, of the stochastic 
differential equation with reflection $(\ref{eq:SDER})$. The 
strong Markov property will be crucial in our argument to 
prove uniqueness of the solution. Our arguments are 
very similar to those of \cite{TW93}, but using the 
constrained martingale problem approach we do not need 
the oscillation estimates $(4.5)$ and $(4.6)$ of \cite{TW93}.

\begin{lemma}\label{th:absorb}
Given a filtered probability space $(\Omega ,{\cal F},\{{\cal F}_
t\},\P )$, a 
standard $\{{\cal F}_t\}$-Brownian motion $W$ and a ${\cal F}_0$-measurable random 
variable $\xi_0$ with compact support in $\overline D-\{0\}$, there exist a 
continuous process $\xi$ and 
a continuous, nondecreasing process $\ell$ such that, setting 
\begin{equation}\theta:=\inf\{t\geq 0:\,\xi (t)=0\},\label{eq:0hitabs}\end{equation}
$\xi$ and $\ell$  satisfy, a.s., 
\begin{eqnarray}
\xi (t)&=&\xi_0+\int_0^tb(\xi (s))ds+\int_0^t\sigma (\xi (s))dW(s
)+\int_0^tg(\xi (s))d\ell (s),\quad 0\leq t<\theta,\non\\
&&\xi (t)\in\overline D-\{0\},\quad\int_0^t{\bf 1}_{\partial D-\{
0\}}(\xi (s))d\ell (s)=\ell (t),\quad 0\leq t<\theta,\label{eq:abs}\end{eqnarray}
and, on the set $\{\theta<\infty \}$, $\xi (t)=0$, $\ell (t)=\ell 
(\theta)$ for $t\geq\theta$. 

Moreover $\xi$ is pathwise unique. 
\end{lemma}

\begin{proof}
Let $\{\delta_k\}$ be a strictly decreasing sequence of positive 
numbers converging to zero, and let $\{D^k\}$ be a sequence 
of bounded domains with ${\cal C}^1$ boundary such that 
$D^k\subset D^{k+1}\subset D$, $\overline {D^k}\cap B_{\delta_k}(
0)^c=\overline D\cap B_{\delta_k}(0)^c$ and $\overline {D^k}\cap 
B_{\delta_k}(0)\subset D^{k+1}.$ 
Also let $g^k:\R^d-\{0\}\rightarrow\R^d$, $k\in\N$, be a locally Lipschitz vector field, of 
unit length on $\partial D^k$, such 
that $g^k(x)=g(x)$ for $x\in\partial D\cap B_{\delta_k}(0)^c$ and that, denoting by 
$n^k(x)$ the unit, inward normal at $x\in\partial D^k$, 
$\inf_{x\in\partial D^k}g^k(x)\cdot n^k(x)>0$. Then we know, by the results 
of \cite{DI93}, that there is a strong solution to 
$(\ref{eq:SDER})$ in $\overline {D^k}$ with direction of reflection $
g^k$ 
and it is pathwise unique. 
Given a Brownian motion $W$, for each random variable $\xi_0$ 
with compact support in $\overline D-\{0\}$, independent of $W$, for all $
k$ 
large enough that $\xi_0$ is supported in $\overline {D^k}$, let $
\xi^k$ be the 
strong solution of $(\ref{eq:SDER})$ in $\overline {D^k}$ with direction 
of reflection $g^k$, with Brownian motion $W$ and 
initial condition $\xi_0$, and $\ell^k$ be the 
corresponding nondecreasing process. 
The sequence of stopping times $\{\theta^k\}$, 
\[\theta^k:=\inf\{t\geq 0:\,\xi^k(t)\in\partial D^k\cap\overline {
B_{\delta_k}(0)}\},\]
is strictly increasing and, setting 
\[\xi (t):=\xi^k(t),\quad\ell (t):=\ell^k(t),\qquad 0\leq t\leq\theta^
k,\]
\[\theta:=\lim_{k\rightarrow\infty}\theta^k,\]
$\xi$ and $\ell$ satisfy $(\ref{eq:abs})$. 

We will now show that, on the set $\{\theta<\infty \}$, 
\begin{equation}\sup_{0\leq t<\theta}\ell (t)<\infty ,\quad
\lim_{t\rightarrow\theta^{-}}\xi (t)=0,\qquad a.s..\label{eq:xiuc}\end{equation}

In fact, let $F$ be the function of Lemma 
\ref{th:bdryfunc}. Then, by Ito's formula, on the set 
$\{\theta<\infty \}$, we have 
\begin{eqnarray*}
\ell (\theta^k)&\leq&\frac 1{c_F}\bigg\{2\|F\|\,+\|AF\|\,\theta^k
+\bigg|\int_0^{\theta^k}\nabla F(\xi (s)^T\sigma (\xi (s))dW(s)\bigg
|\bigg\}\\
\ &=&\frac 1{c_F}\bigg\{2\|F\|\,+\|AF\|\theta^k+\bigg|\int_0^{\theta^
k}{\bf 1}_{\{s<\theta\}}\nabla F(\xi (s)^T\sigma (\xi (s))dW(s)\bigg
|\bigg\}.\end{eqnarray*}
Since the process ${\bf 1}_{\{s<\theta\}}\nabla F(\xi (s)^T\sigma 
(\xi (s))$ is predictable and 
bounded, the process\hfill\break
 $\int_0^t{\bf 1}_{\{s<\theta\}}\nabla F(\xi (s)^T\sigma (\xi (
s))dW(s)$ is a.s. 
continuous. Therefore, on the set $\{\theta<\infty \}$, the limit 
$\lim_{k\rightarrow\infty}\int_0^{\theta^k}{\bf 1}_{\{s<\theta\}}
\nabla F(\xi (s)^T\sigma (\xi (s))dW(s)$ exists and is finite 
a.s.. This yields the first assertion in $(\ref{eq:xiuc})$ and 
allows to define, on the set $\{\theta<\infty \}$, $\ell (\theta):=\sup_{0\leq t<\theta}\ell (t).$
The second assertion follows by observing that, on the set 
$\{\theta<\infty \}$, all the terms on the r.h.s. of the first equation 
in $(\ref{eq:abs})$ are uniformly continuous 
on $[0,\theta)$. 
\end{proof}

\begin{theorem}\label{th:exist} 
Under Conditions \ref{D}, \ref{G} and \ref{A}, 
for each $\nu\in {\cal P}(\overline D)$, there exists a strong Markov, natural solution of the 
constrained 
martingale problem for $(A,D,B,\Xi )$, with initial distribution $
\nu$. 
\end{theorem}

\begin{proof}
We will first construct a solution of the controlled 
martingale problem. 
Let $g^0$ be a unit vector in $G(0)\cap {\cal K}$ (Condition \ref{G} (iii)). 
Then, by Condition 
\ref{D}(i), for $\rho$ small enough $\rho g^0\in G(0)\cap D$. 
Let $\{\rho_n\}$ be 
a decreasing sequence of positive numbers such that 
$\rho_n\rightarrow 0$ and $\rho_ng^0\in G(0)\cap D$ for all $n$.  
Consider the stochastic differential equation with 
reflection 
\begin{eqnarray}
X^n(t)&=&X^n_0+\hskip-.01in\int_0^tb(X^n(s))ds+\hskip-.01in\int_0^t\sigma (X^n(s))dW(s)+\hskip-.01in\int_
0^tg(X^n(s))dl^n(s)+\rho_ng^0L^n(t),\non\\
&&X^n(t)\in\overline D-\{0\},\quad l^n\mbox{\rm \ non decreasing}
,\quad\int_0^t{\bf 1}_{\partial D-\{0\}}(X^n(s))dl^n(s)=l^n(t),\quad\label{
eq:sSDER}\\
&&\qquad\qquad L^n(t)=\#\{s\leq t:X^n(s^{-})=0\},\non\end{eqnarray}
where $\#$ denotes cardinality, $W$ is a standard Brownian motion, 
$X^n_0$ is a random variable independent of $W$, with law $\nu_n$ 
defined as 
\[\nu_n(C):=\nu (\{0\})\delta_{x^n}(C)+\nu (\overline D-\{0\})\frac {
\nu (C\cap K_n)}{\nu (K_n)},\qquad\forall C\in {\cal B}(\overline 
D),\]
where $\{x^n\}$ is a sequence of points in $\overline D-\{0\}$ converging 
to $0$ and $\{K_n\}$ is an increasing sequence of compact sets of $\R^
d$ such that 
$\bigcup_{n=1}^{\infty}K_n=\overline D-\{0\}$. Then $\nu_n$ has compact support in $\overline 
D-\{0\}$ 
and $\{\nu_n\}$ converges weakly to $\nu$. 
Existence of $X^n$ follows from Lemma \ref{th:absorb}. 
 
Define 
\[\lambda^n_0(t):=\inf\{s:s+l^n(s)+\rho^nL^n(s)>t\},\]
\[\Lambda_1^n([0,t]\times C):=\int_0^t{\bf 1}_C(g(X^n(\lambda^n_0
(s)))\,dl^n(\lambda_0^n(s)))+{\bf 1}_C(g^0)\rho_nL^n(\lambda_0^n(
s)),\]
\[B^nf(x,u):=u\cdot\nabla f(x){\bf 1}_{\partial D-\{0\}}(x)+(\rho_
n)^{-1}[f(x+\rho_nu)-f(x)]{\bf 1}_{\{0\}}(x),\]
and  $Y^n(t):=X^n(\lambda_0^n(t))$, in particular $Y^n(0)=X^n_0$.  Then for 
each $f\in {\cal C}^2(\overline D)$,
\[f(Y^n(t))-f(Y^n(0))-\int_0^tAf(Y^n(s))d\lambda_0^n(s)-\int_{[0,
t]\times U}B^nf(Y^n(s^{-}),u)\Lambda^n_1(ds\times du)\]
is a martingale with respect to $\{{\cal F}^{Y^n,\lambda_0^n,\Lambda_
1^n}_t\}$. 
$(Y^n,\lambda_0^n,\Lambda_1^n)$ is not a solution of the controlled martingale 
problem for $(A,D,B^n,\Xi )$ in the precise sense of 
\cite{CK19} and \cite{Kur91}, because $B^nf$ is not 
continuous on $\Xi$ and because we can only say
\[t\leq\lambda_0^n(t)+\Lambda_1^n([0,t]\times U)\leq t+\rho_n.\]
However, the same relative 
compactness arguments as, for example, in Lemma 
2.8 of \cite{CK19} apply and, since the law of $Y^n(0)$, $\nu_n$,  
converges to $\nu$, any limit 
point of $\{(Y^n,\lambda_0^n,\Lambda_1^n)\}$ will be a solution of the controlled 
martingale problem for $(A,D,B,\Xi )$ with initial distribution 
$\nu$.  

Then the assertion follows from Lemma 
\ref{th:bdryfunc}, Lemma \ref{th:takeoff}, Lemma 3.1 of 
\cite{CK19} and Corollary 4.12 a) of \cite{CK19}.

\end{proof}

\subsection{Uniqueness}\label{sectionuniq}

We state preliminarly the following lemma, which holds 
for constrained martingale problems in general.

\begin{lemma}\label{th:unq-by-1dim}
Suppose that for each $\nu\in {\cal P}(\overline D)$ any two strong Markov, natural solutions of the 
constrained martingale problem for $(A,D,B,\Xi )$ with initial distribution 
$\nu$ have the same one-dimensional distributions. Then, for each $
\nu\in {\cal P}(\overline D)$, 
any two strong Markov, 
natural solutions of the constrained martingale problem for $(A,D
,B,\Xi )$ with initial distribution 
$\nu$ have the same distribution. 
\end{lemma}

\begin{proof}
The proof of point (a) of Theorem 4.4.2 of \cite{EK86} 
carries over. 
\end{proof}
\medskip

Our uniqueness result is formulated in the following 
theorem. By Theorem \ref{th:equiv}, it is equivalent 
to uniqueness in distribution for the solution of 
$(\ref{eq:SDER})$. 

\begin{theorem}\label{th:unique}
Under Conditions \ref{D}, \ref{G} and \ref{A}, for every  
$\nu\in {\cal P}(\overline D)$, there is a unique natural 
solution of the constrained martingale problem 
for $(A,D,B,\Xi )$ with initial distribution 
$\nu$. The solution is a strong Markov process. 
\end{theorem}

The proof of Theorem \ref{th:unique} requires several 
steps. The first building block is the following lemma.
For $\delta^{*}$, $0<\delta^{*}\leq\bar{\delta}$, where $\bar{\delta}$ is as in Lemma \ref{th:extime}, let 
\begin{equation}D_n:=\overline D\cap B_{\delta^{*}2^{-2n}},\quad 
E_n:=\{x\in\overline D:|x|=\delta^{*}2^{-2n}\},\quad n\geq 0.\label{eq:Dn}\end{equation}
Let $X$ be a solution of 
the constrained martingale problem for $(A,D,B,\Xi )$. Define 
recursively 
\begin{eqnarray}
\tau_0^n:=0,\quad\vartheta^n_0:=\vartheta:=\inf\{t\geq 0:X(t)=0
\},\quad\tau^n:=\inf\{t\geq 0:X(t)\in E_n\},\;n\geq 0,\non\\
\hfill\\
\tau_l^n:=\inf\{t\geq\vartheta_{l-1}^n:X(t)\in E_n\},\quad\vartheta_
l^n:=\inf\{t\geq\tau_l^n:X(t)=0\},\qquad l\geq 1,\;n\geq 0.\non\label{eq:ballcycles}
\end{eqnarray}

\begin{lemma}\label{th:unq-by-hit}
Suppose that, for some $\delta^{*}$, $0<\delta^{*}\leq\bar{\delta}$, the hitting distributions $
\{\mu_n\}$ defined by 
\begin{equation}\mu_n(C):=\mathbb P\{X(\tau^n)\in C\},\quad C\in 
{\cal B}(E_n),\quad n\geq 0,\label{eq:hitdist}\end{equation}
are the same for any strong Markov, natural solution $X$ 
of the constrained martingale problem for $(A,D,B,\Xi )$ 
starting at $0$. 

Then, for each $\nu\in {\cal P}(\overline D)$,  there is a unique natural solution of 
the constrained martingale problem for $(A,D,B,\Xi )$ 
with initial distribution $\nu$. 
\end{lemma}

\begin{proof}
Lemma \ref{th:takeoff} allows to apply Corollary 4.13 of 
\cite{CK19}. Therefore it is enough to prove 
uniqueness among strong Markov, natural solutions of 
the constrained martingale problem for $(A,D,B,\Xi )$. 
Let $X$ be such a solution with initial distribution $\nu$.  

For $\eta >0$ and $f\in {\cal C}(\overline D)$ that vanishes in a neighborhood of 
the origin, define, for each $n\geq 0$, 
\[R^n_{\eta}f:=\E\Big[\int_0^{^{\vartheta}}e^{-\eta t}f(X(t))dt\Big
]+\E\Big[\sum_{l=1}^{\infty}\prod_{m=0}^{l-1}e^{-\eta (\vartheta_
m^n-\tau_m^n)}\int_{\tau_l^n}^{\vartheta_l^n}e^{-\eta (t-\tau_l^n
)}f(X(t))dt\Big],\]
with the convention $e^{-\infty}=0$.
The hypotheses ensure that the distribution of $X(\tau_l^n)$ is 
$\mu_n$ for all $l$. Then, by the strong Markov property, 
Lemma \ref{th:absorb} and Theorem \ref{th:equiv}, for 
each $l$ the factors 
of the product in the infinite sum are independent, with 
distributions  determined by the initial distribution $\nu$ 
and by the unique distribution of 
$X^n(\cdot\wedge\theta^n)$, where $X^n$ is a solution of the constrained 
martingale problem with initial distribution $\mu_n$ and 
$\theta^n$ is the first time $X^n$ hits zero.  Consequently, each term on 
the right hand side is uniquely determined. 
The independence and the fact that $\{\vartheta_m^n-\tau_m^n\}_{m
\geq 1}$ are 
identically distributed 
imply also that the series is convergent and 
has finite expectation.  

Let ${\cal T}^n:=[0,\vartheta)\cup\cup_{l=1}^{\infty}[\tau_l^n,
\vartheta_l^n)$.  Then 
$R_{\eta}^nf$ can be written as 
\[R_{\eta}^nf=\mathbb E\left[\int_0^{\infty}{\bf 1}_{{\cal T}^n}(
t)e^{-\eta\int_0^t{\bf 1}_{{\cal T}^n}(s)ds}f(X(t))dt\right].\]
We are going to show that 
\[R_{\eta}^nf\;\rightarrow_{n\rightarrow\infty}\mathbb E\left[\int_
0^{\infty}e^{-\eta t}f(X(t))dt\right]:=R_{\eta}f.\]
First of all note that, for $n$ large enough, depending only on $
f$, 
\[\int_0^{\infty}{\bf 1}_{{\cal T}^n}(t)\,e^{-\eta\int_0^t{\bf 1}_{
{\cal T}^n}(s)ds}f(X(t))dt=\int_0^{\infty}e^{-\eta\int_0^t{\bf 1}_{
{\cal T}^n}(s)ds}f(X(t))dt,\quad a.s.,\]
\[\int_0^{\infty}{\bf 1}_{{\cal T}^n}(t)\,e^{-\eta\int_0^t{\bf 1}_{
{\cal T}^n}(s)ds}|f(X(t))|dt=\int_0^{\infty}e^{-\eta\int_0^t{\bf 1}_{
{\cal T}^n}(s)ds}|f(X(t))|dt,\quad a.s..\]

By Lemma \ref{th:occtime}, for each $t\geq 0$, 
\begin{equation}e^{-\eta\int_0^t{\bf 1}_{{\cal T}^n}(s)ds}\rightarrow_{
n\rightarrow\infty}e^{-\eta t},\qquad a.s.,\label{eq:occtlim}\end{equation}
therefore we only need to show that we can pass to the 
limit in the time integral and under the expectation. 
Now we have 
\[\int_0^{\infty}{\bf 1}_{{\cal T}^n}(s)ds=\vartheta+\sum_{l=1}^{
\infty}(\vartheta_l^n-\tau_l^n)=\infty ,\quad a.s.,\;\forall n,\]
because the random variables in the infinite sum in the right 
hand side are positive i.i.d. random variables. Hence, by $(\ref{eq:occtlim})$, 
\begin{eqnarray*}
&&\lim_{t_0\rightarrow\infty}\limsup_{n\rightarrow\infty}\int_{t_
0}^{\infty}e^{-\eta\int_0^t{\bf 1}_{{\cal T}^n}(s)ds}|f(X(t))|dt\non\\
&&=\lim_{t_0\rightarrow\infty}\limsup_{n\rightarrow\infty}\int_{t_
0}^{\infty}{\bf 1}_{{\cal T}^n}(t)\,e^{-\eta\int_0^t{\bf 1}_{{\cal T}^
n}(s)ds}|f(X(t))|dt\\
&&\leq\lim_{t_0\rightarrow\infty}\limsup_{n\rightarrow\infty}\frac {
\|f\|}{\eta}e^{-\eta\int_0^{t_0}{\bf 1}_{{\cal T}^n}(s)ds}=0\qquad 
a.s.,\non\label{eq:unint}\end{eqnarray*}
so that, almost surely, the sequence $\big\{e^{-\eta\int_0^t{\bf 1}_{
{\cal T}^n}(s)ds}f(X(t))\big\}$ 
is uniformly integrable in time and we can pass to the 
limit in the time integral. Analogously 
\[\int_0^{\infty}e^{-\eta\int_0^t{\bf 1}_{{\cal T}^n}(s)ds}|f(X(t
))|dt\leq\frac {\|f\|}{\eta},\quad a.s.,\]
and we can also pass to the limit under the expectation. 

The class of continuous functions on $\overline D$ that vanish 
in a neighborhood of the origin is separating for the 
probability measures on $\overline D$, therefore 
it follows that the one-dimensional distributions of $X$ are 
uniquely determined, and hence the assertion, by Lemma 
\ref{th:unq-by-1dim}. 
\end{proof}
\medskip

The next lemma shows how Theorem \ref{th:ergodic} 
comes into play in verifying the assumption of Lemma 
\ref{th:unq-by-hit}. Recall that, by Lemma 
\ref{th:absorb}, for all $n\geq 0$, $x\in\overline D-\{0\}$, $f\in 
{\cal C}(E_n)$, the expectation 
\[\E[f(X(\tau^n))\mathbf{1}_{\tau^n<\vartheta}|X(0)=x]\]
is the same for all natural solutions $X$ of the constrained martingale problem 
for $(A,D,B,\Xi )$. 

\begin{lemma}\label{th:hit-by-erg}
For $k\geq 1$, let $Q_k$ be the subprobability transition operator defined by 
\begin{equation}Q_kf(x):=\E[f(X(\tau^{k-1}))\mathbf{1}_{\tau^{k-1}
<\vartheta}|X(0)=x],\quad x\in E_k,\quad f\in {\cal C}(E_{k-1})
.\label{eq:Qk}\end{equation}
Then, for every  strong Markov, natural solution $X$ 
of the constrained martingale problem for $(A,D,B,\Xi )$ 
starting at $0$,  for every $n\geq 0$, 
\[\E[f(X(\tau^n))]=\frac {\int Q_{n+k}\cdots Q_{n+1}f(x)\mu_{n+k}
(dx)}{\int Q_{n+k}\cdots Q_{n+1}1(x)\mu_{n+k}(dx)},\quad\forall f
\in {\cal C}(E_n),\quad\forall k\geq 1,\]
where $\mu_{n+k}$ is defined by $(\ref{eq:hitdist})$. 
\end{lemma}

\begin{proof}
Note that, since $X$ starts at $0$, $\tau^n=\tau_1^n$ for all $n\geq 
0$. By the strong 
Markov property we have, for $n\geq 0$,   
\begin{eqnarray*}
&\E&[f(X(\tau^n))]\\
&=&\E[f(X(\tau^n))\mathbf{1}_{\vartheta_1^{n+k}<\tau^{n+k-1}}]\\
&&+\E[f(X(\tau^n))\mathbf{1}_{\vartheta_1^{n+k}>\tau^{n+k-1}}\mathbf{
1}_{\vartheta_1^{n+k-1}<\tau^{n+k-2}}]\\
&&+\;\cdots\cdots\cdots\cdots\cdots\cdots\\
&&+\E[f(X(\tau^n))\mathbf{1}_{\vartheta_1^{n+k}>\tau^{n+k-1}}\cdots\mathbf{
1}_{\vartheta_1^{n+2}>\tau^{n+1}}\mathbf{1}_{\vartheta_1^{n+1}<\tau^
n}]\\
&&+\E[f(X(\tau^n))\mathbf{1}_{\vartheta_1^{n+k}>\tau^{n+k-1}}\cdots\mathbf{
1}_{\vartheta_1^{n+2}>\tau^{n+1}}\mathbf{1}_{\vartheta_1^{n+1}>\tau^
n}]\\
\\
&=&\E[\E[f(X(\tau^n))|{\cal F}_{\vartheta_1^{n+k}}]\mathbf{1}_{\vartheta_
1^{n+k}<\tau^{n+k-1}}]\\
&&+\E[\E[f(X(\tau^n))|{\cal F}_{\vartheta_1^{n+k-1}}]\mathbf{1}_{
\vartheta_1^{n+k}>\tau^{n+k-1}}\mathbf{1}_{\vartheta_1^{n+k-1}<\tau^{
n+k-2}}]\\
&&+\;\cdots\cdots\cdots\cdots\cdots\cdots\\
&&+\E[\E[f(X(\tau^n))|{\cal F}_{\vartheta_1^{n+1}}]\mathbf{1}_{\vartheta_
1^{n+k}>\tau^{n+k-1}}\cdots\mathbf{1}_{\vartheta_1^{n+2}>\tau^{n+
1}}\mathbf{1}_{\vartheta_1^{n+1}<\tau^n}]\\
&&+\E[\mathbf{1}_{\vartheta_1^{n+k}>\tau^{n+k-1}}\cdots\mathbf{1}_{
\vartheta_1^{n+2}>\tau^{n+1}}\E[f(X(\tau^n))\mathbf{1}_{\vartheta_
1^{n+1}>\tau^n}|{\cal F}_{\tau^{n+1}}]]\\
&=&\E[f(X(\tau^n))]\big(1-\E[\E[\mathbf{1}_{\vartheta_1^{n+k}>\tau^{
n+k-1}}\cdots\E[\mathbf{1}_{\vartheta_1^{n+1}>\tau^n}|{\cal F}_{\tau^{
n+1}}]\cdots |{\cal F}_{\tau^{n+k}}]]\big)\\
&&+\E[\E[\mathbf{1}_{\vartheta_1^{n+k}>\tau^{n+k-1}}\cdots\E[\mathbf{
1}_{\vartheta_1^{n+2}>\tau^{n+1}}\E[f(X(\tau^n))\mathbf{1}_{\vartheta_
1^{n+1}>\tau^n}|{\cal F}_{\tau^{n+1}}]|{\cal F}_{\tau^{n+2}}]\cdots 
|{\cal F}_{\tau^{n+k}}]]\\
&=&\E[f(X(\tau^n))]\big(1-\E[\E[\mathbf{1}_{\vartheta_1^{n+k}>\tau^{
n+k-1}}\cdots\E[\mathbf{1}_{\vartheta_1^{n+1}>\tau^n}|X(\tau^{n+1}
)]\cdots |X(\tau^{n+k})]]\big)\\
&&+\E[\E[\mathbf{1}_{\vartheta_1^{n+k}>\tau^{n+k-1}}\cdots\E[\mathbf{
1}_{\vartheta_1^{n+2}>\tau^{n+1}}\E[f(X(\tau^n))\mathbf{1}_{\vartheta_
1^{n+1}>\tau^n}|X(\tau^{n+1})]\cdots |X(\tau^{n+k})]]\\
&=&\E[f(X(\tau^n))]\big(1-\E[Q_{n+k}\cdots Q_{n+1}1(X(\tau^{n+k})
]\big)\\
&&+\E[Q_{n+k}\cdots Q_{n+1}f(X(\tau^{n+k})].\end{eqnarray*}
\end{proof}\

\noindent {\bf Proof of Theorem \ref{th:unique}}
Suppose that, for some $\delta^{*}$, $0<\delta^{*}\leq\bar{\delta}$ and for each fixed $
n\geq 0$, the sequence 
of subprobability transition functions $\{Q_{n+k}\}_{k\geq 1}$ defined in 
Lemma \ref{th:hit-by-erg} satisfies 
the assumptions of Theorem \ref{th:ergodic}. Then, for every 
$f\in {\cal C}(E_n)$ there exists a constant $C_n(f)$ such that 
\[\E[f(X(\tau^n))]=C_n(f)\]
for any strong Markov, natural solution $X$ of the constrained martingale problem 
for $(A,D,B,\Xi )$ starting at $0$, and the theorem follows 
from Lemma \ref{th:unq-by-hit}. Therefore we are 
reduced to verifying that 
the sequence of subprobability transition functions $\{Q_{n+k}\}_{
k\geq 1}$ satisfies 
the assumptions of Theorem \ref{th:ergodic}. This is the 
object of the next two subsections. More precisely, in 
Section \ref{sectionhittimes} we 
show that the subprobability transition functions $\{Q_{n+k}\}_{k
\geq 1}$
are not identically zero and that condition (i) is verified and in Section 
\ref{sectionhitdist} we show that condition (ii) is verified. 
\hfill $\Box$\medskip

\subsection{Estimates on hitting times}\label{sectionhittimes}

In this subsection we prove that, for $\delta^{*}$ small enough, 
for each $n\geq 0$, the subtransition functions $\{Q_{n+k}\}_{k\geq 
1}$ defined 
by $(\ref{eq:Qk})$, $(\ref{eq:ballcycles})$ and $(\ref{eq:Dn})$ 
are not identically zero and 
verify assumption (i) of Theorem 
\ref{th:ergodic}. 

We have, for $x\in E_{n+k}$, 
\begin{eqnarray}
Q_{n+k}1(x)=\P (\tau^{n+k-1}<\vartheta|X(0)=x),\quad
\non\\
\quad\hfill \\
\qquad Q_{n+k}\cdots Q_{n+1}1(x)=\P (\tau^n<\vartheta|X(0)=x),\non\label{eq:i-to-hitprob}\end{eqnarray}
where $X$ is a strong Markov, natural solution of 
the constrained martingale problem for $(A,D,B,\Xi )$ and the 
right hand sides in the above equalities are uniquely determined by Lemma 
\ref{th:absorb}. Recall that we use the notation 
\[\P^x(\cdot ):=\P (\cdot |X(0)=x),\quad\E^x[\cdot ]:=\E[\cdot |X
(0)=x].\]
The fact that $Q_{n+k}$ is not identically zero holds if, more 
generally, for $0<\delta\leq\delta^{*}$,
\begin{equation}\P^x(\tau^{\delta}<\vartheta)>0,\quad\forall x\in\overline 
D,\,0<|x|<\delta,\label{eq:non0hitprob}\end{equation}
where $\tau^{\delta}:=\tau^{X,\delta}$ is defined by $(\ref{eq:taudelta})$. 
Analogously, supposing $\P^x(\tau^{\delta}<\vartheta)>0$ for all 
$x\in\overline D,\,0<|x|<\delta$, $0<\delta\leq\delta^{*}$, assumption (i) of Theorem 
\ref{th:ergodic} is verified if there exists 
$c_0>0$ such that, for all $0<\delta\leq\delta^{*}$, 
\begin{equation}\inf_{x,y\in \overline{D}:\,0<|x|=|y|<\delta}\frac {\P^x(\tau^{\delta}
<\vartheta)}{\P^y(\tau^{\delta}<\vartheta)}\geq c_0,\label{eq:hitprob}\end{equation}

The proof of $(\ref{eq:non0hitprob})$ and $(\ref{eq:hitprob})$ 
is based on estimating $\P^x(\tau^{\delta}<\vartheta)$ by means of 
suitable auxiliary functions (Lemmas \ref{th:hitprob} and \ref{th:auxfunc}). 
These auxiliary functions are constructed by elaborating 
on some functions introduced by \cite{VW85} and 
\cite{KW91} in the study of the reflecting Brownian 
motion in a cone with radially constant direction of 
reflection. (See 
Appendices  \ref{appcone} and \ref{appaux}.) 

Let $\alpha^{*}$ be defined as in Theorems \ref{th:funcd2} and 
\ref{th:funcd3} for ${\cal K}$ and $\bar {g}$ given by Conditions \ref{D} 
and \ref{G}. 

\begin{lemma}\label{th:auxfunc}
There exists $\delta^{*}>0$ such that:
\begin{itemize}
\item[(i)]If $\alpha^{*}\leq 0$, there exists a function $V\in {\cal C}^
2(\overline D-\{0\})$ such 
that 
\begin{equation}\lim_{x\in\overline D,\,x\rightarrow 0}V(x)=\infty 
,\label{eq:V-zero}\end{equation}
\begin{equation}\nabla V(x)\cdot g(x)\leq 0,\quad\forall x\in\big
(\partial D-\{0\}\big)\cap\overline {B_{\delta^{*}}(0)}\label{eq:nablaV-}\end{equation}
\begin{equation}AV(x)\leq 0,\quad\forall x\in\big(\overline D-\{0
\}\big)\cap\overline {B_{\delta^{*}}(0)}.\label{eq:AV-}\end{equation}
\item[(ii)]If $0<\alpha^{*}<1$, there exist two functions 
$V_1,V_2\in {\cal C}^2(\overline D-\{0\})$ such that 
\begin{equation}\begin{array}{c}
V_1(x)>0,\quad V_2(x)>0,\quad\mbox{\rm for }x\in\big(\overline D-
\{0\}\big)\cap\overline {B_{\delta^{*}}(0)},\\
\\
\lim_{x\in\overline D,\,x\rightarrow 0}V_1(x)=\lim_{x\in\overline 
D,\,x\rightarrow 0}V_2(x)=0,\\
\\
\quad\inf_{0<\delta\leq\delta^{*}}\frac {\inf_{|x|=\delta}V_1(x)}{\sup_{
|x|=\delta}V_2(x)}>0,\quad\inf_{0<\delta\leq\delta^{*}}\frac {\inf_{
|x|=\delta}V_2(x)}{\sup_{|x|=\delta}V_1(x)},>0\end{array}
\label{eq:V+zero}\end{equation}
\begin{equation}\nabla V_1(x)\cdot g(x)\geq 0,\quad\nabla V_2(x)\cdot 
g(x)\leq 0,\quad\forall x\in\big(\partial D-\{0\}\big)\cap\overline {
B_{\delta^{*}}(0)}\label{eq:nablaV+}\end{equation}
\begin{equation}AV_1(x)\geq 0,\quad AV_2(x)\leq 0,\quad\forall x\in\big
(\overline D-\{0\}\big)\cap\overline {B_{\delta^{*}}(0)}.\label{eq:AV+}\end{equation}
\end{itemize}
\end{lemma}

\begin{proof}
See Appendix \ref{appaux}.
\end{proof}

\begin{lemma}\label{th:hitprob}
Assume Conditions \ref{D}, \ref{G} and \ref{A}. Then $\alpha^{*}<1$. 

For a natural solution, $X$, of the constrained 
martingale problem for $(A,D,B,\Xi )$, let $\vartheta$ be defined by 
$(\ref{eq:ballcycles})$ and $\tau^{\delta}$ be defined by $(\ref{eq:taudelta})$. 
Then there exists $\delta^{*}>0$ such that, for every $x\in\overline D$, 
$0<|x|<\delta\leq\delta^{*}$, $P^x\big(\tau^{\delta}<\infty\big)=1$ and: 
\begin{itemize}
\item[(i)]if $\alpha^{*}\leq 0$, 
\[\P^x\big(\tau^{\delta}<\vartheta\big)=1.\]
\item[(ii)]if $0<\alpha^{*}<1$, 
\[\P^x(\tau^{\delta}<\vartheta)>0.\]
Moreover there exists $c_0>0$ such that, for all $0<\delta\leq\delta^{
*}$,  
\[\inf_{x,y\in \overline{D}:\,0<|x|=|y|<\delta}\frac {\P^x(\tau^{\delta}<\vartheta)}
{\P^y(\tau^{\delta}<\vartheta)}\geq c_0.\]
\end{itemize}
\end{lemma}

\begin{proof}
$\alpha^{*}<1$ by Theorem \ref{th:alpha}. 
Let $\delta^{*}$, $V$, $V_1$, $V_2$ be as in Lemma \ref{th:auxfunc}.

Let $\alpha^{*}\leq 0$. Of course we can always suppose that $V$ is 
nonnegative on $\big(\overline D-\{0\}\big)\cap\overline {B_{\delta^{
*}}(0)}$. 
By applying Ito's formula to the function $V$, we 
have, for $\delta\leq\delta^{*}$, for every fixed $x\in\big(\overline 
D-\{0\}\big)\cap B_{\delta}(0)$ and $\epsilon <|x|$, 
\[\E^x[V(X(\tau^{\epsilon}\wedge\tau^{\delta})]\leq V(x),\]
which yields 
\[\inf_{|y|=\epsilon}V(y)\P^x\big(\tau^{\epsilon}<\tau^{\delta}\big
)\leq V(x),\]
and hence, by letting $\epsilon\rightarrow 0$, by $(\ref{eq:V-zero})$, 
\[\P^x\big(\vartheta<\tau^{\delta}\big)=0.\]

If $0<\alpha^{*}<1$, by applying Ito's formula to $V_1$, we obtain, 
for $x\in\big(\overline D-\{0\}\big)\cap B_{\delta}(0)$ and $\epsilon 
<|x|$, 
\[\E^x[V_1(X(\tau^{\delta}\wedge\tau^{\epsilon}))]\geq V_1(x),\]
which yields 
\[\sup_{|u|=\delta}V_1(u)\,\P^x\big(\tau^{\delta}<\tau^{\epsilon}\big
)+\sup_{|u|=\epsilon}V_1(u)\P^x\big(\tau^{\epsilon}<\tau^{\delta}\big
)\geq V_1(x),\]
and hence, by letting $\epsilon\rightarrow 0$, 
\begin{equation}\sup_{|u|=\delta}V_1(u)\,\P^x\big(\tau^{\delta}<\vartheta
\big)\geq V_1(x),\label{eq:alpha+lb}\end{equation}
which yields the first assertion in part (ii) of the 
thesis. Analogously, by applying Ito's formula to $V_2$ we get, 
\[\inf_{|u|=\delta}V_2(u)\,\P^x\big(\tau^{\delta}<\tau^{\epsilon}\big
)+\inf_{|u|=\epsilon}V_2(u)\,\P^x\big(\tau^{\epsilon}<\tau^{\delta}\big
)\leq V_2(x),\]
and hence, by letting $\epsilon\rightarrow 0$, 
\begin{equation}\inf_{|u|=\delta}V_2(u)\,\P^x\big(\tau^{\delta}<\vartheta
\big)\leq V_2(x).\label{eq:alpha+ub}\end{equation}
Combining $(\ref{eq:alpha+lb})$ and $(\ref{eq:alpha+ub})$, we get, for 
$x,y\in\big(\overline D-\{0\}\big)\cap B_{\delta}(0)$ with $|x|=|
y|$, 
\[\frac {\P^x\big(\tau^{\delta}<\vartheta\big)}{\P^y\big(\tau^{
\delta}<\vartheta\big)}\geq\frac {V_1(x)}{V_2(y)}\,\frac {\inf_{
|u|=\delta}V_2(u)}{\sup_{|u|=\delta}V_1(u)}\geq\inf_{0<\delta\leq
\delta^{*}}\frac {\inf_{|u|=\delta}V_1(u)}{\sup_{|u|=\delta}V_2(u
)}\inf_{0<\delta\leq\delta^{*}}\frac {\inf_{|u|=\delta}V_2(u)}{\sup_{
|u|=\delta}V_1(u)}\,>0.\]
\end{proof}

\subsection{Estimates on hitting 
distributions}\label{sectionhitdist}

In this subsection we verify assumption (ii) of 
Theorem \ref{th:ergodic} (Lemma \ref{th:coupling}). 
Lemma \ref{th:coupling} follows essentially from 
the fact that, for any $x,\,\tilde {x}\in E_n$, we can construct, on 
the same probability space, two strong Markov, natural solutions of the 
constrained martingale problem for $(A,D,B,\Xi )$, starting at 
$x$ and $\tilde {x}$, such that the probability that they hit $E_{
n-1}$ 
before the origin and that 
they couple before hitting $E_{n-1}$ (i.e. 
that their paths agree, up to a time shift, for some time 
before they hit $E_{n-1}$) is larger than some $\epsilon_0>0$ 
independent of $x$ and $\tilde {x}$ and of $n$. 
The construction is based on a 
coupling result of \cite{CK18} and on a uniform lower bound on 
the probability that a natural solution 
of the constrained martingale problem for $(A,D,B,\Xi )$ 
starting on $E_n$ hits the intermediate layer 
$\{x\in\overline D;\,|x|=2^{-2n+1}\delta^{*}\}$ in the set 
${\cal O}^n:=\{x\in D:\,2^{2n}x\in {\cal O}\}$, where ${\cal O}$ is an arbitrary open set 
such that ${\cal O}\cap {\cal K}\cap\partial B_{2\delta^{*}}(0)\neq
\emptyset$. In turn this 
uniform lower bound is proved 
by showing that, for any natural 
solution of the constrained martingale problem for 
$(A,D,B,\Xi )$, $X$, the rescaled process $2^{2n}\,X(2^{-4n}\cdot 
)$ 
converges to the reflecting Brownian motion in $\overline {{\cal K}}$ with 
direction of reflection $\bar {g}$ (Lemma \ref{th:scaling}) and 
by the support theorem 
of \cite{KW91}.  
Existence and uniqueness of the reflecting Brownian motion in $\overline {
{\cal K}}$ with 
direction of reflection $\bar {g}$ has been proved in \cite{VW85} and \cite{KW91}, assuming only Condition 
\ref{D} and Conditions \ref{G} (i) and (ii). We 
show in Appendix \ref{appcone} (Theorem 
\ref{th:cone-exist}) that, if Conditions \ref{G} 
(iii) and (iv) are verified, this reflecting Brownian motion 
is the unique natural solution of the constrained 
martingale problem for $(\frac 12\Delta ,{\cal K},B,\bar{\Xi })$, where $
\Delta$ is the 
Laplacian operator, ${\cal D}(\Delta ):={\cal C}^2_b(\overline {{\cal K}}
)$, and 
\begin{eqnarray}\label{eq:cone-cmp}
&\bar{\Xi }:=\{(x,u)\in\partial {\cal K}\times S^{d-1}:u\in\bar {
G}(x)\},&\non\\
&\bar {G}(x):=\{\eta\bar {g}(x),\,\eta\geq 0\},\mbox{\rm \ for }x
\in\partial {\cal K}-\{0\},\quad\bar {G}(0):=G(0),&\\
&B:{\cal C}^2_b(\overline {{\cal K}})\,\rightarrow\,{\cal C}(\bar{
\Xi }),\quad Bf(x,u):=\nabla f\cdot u.&\non\end{eqnarray}
In particular the reflecting Brownian motion is a 
semimartingale (see Remark \ref{re:semimg}.) The reflecting Brownian motion in $\overline {
{\cal K}}$ with 
direction of reflection $\bar {g}$ will be denoted by $\bar {X}$. 

Recall that we are assuming 
Conditions \ref{D}, \ref{G} and \ref{A} and that, for a 
stochastic process $Z$, we use the notation 
\[\P^z(\cdot ):=\P (\cdot |Z(0)=z),\quad\E^z[\cdot ]:=\E[\cdot |Z
(0)=z].\]
Let $X$ be a natural solution of the constrained 
martingale problem for $(A,D,B,\Xi )$. Define,  
\begin{equation}\sigma^{n-1}:=\inf\{t\geq 0:|X(t)|=2^{-2n+1}\delta^{
*}\},\quad\quad n\geq 1,\label{eq:sigman}\end{equation}
and note that $\sigma^{n-1}$ is the hitting time of the surface 
``halfway'' between $E_n$ and $E_{n-1}$. Recall that $\vartheta$ is defined 
by $(\ref{eq:ballcycles})$.

\begin{lemma}\label{th:scaling}
For any sequence $\{x^n\}\subseteq\overline D-\{0\}$ such that $\{
2^{2n}x^n\}$ converges to 
some $\bar {x}\in\overline {{\cal K}}-\{0\}$, let $X^{x^n}$ be a natural solution of the constrained 
martingale problem for $(A,D,B,\Xi )$ starting at $x^n$ and 
$\bar {X}^{\bar {x}}$ be the reflecting Brownian motion in $\overline {
{\cal K}}$ with direction 
of reflection $\bar {g}$, starting at $\bar {x}$. Then 
\begin{equation}2^{2n}\,X^{x^n}(2^{-4n}\cdot )\,\stackrel {{\cal L}}{
\rightarrow}\bar {X}^{\bar {x}}(\cdot ).\label{eq:converg}\end{equation}
In particular, for any open set ${\cal O}$ such that 
${\cal O}\cap {\cal K}\cap\partial B_{2\delta^{*}}(0)\neq\emptyset$, there exists $
\eta_0=\eta_0({\cal O})>0$ such that, 
for $|x^n|=2^{-2n}\delta^{*}$, $\{2^{2n}x^n\}$ converging to $\bar {
x}$, and ${\cal O}^n:=\{x:\,2^{2n}x\in {\cal O}\}$, 
\begin{equation}\liminf_{n\rightarrow\infty}\P^{x^n}(\sigma^{n-1}
<\vartheta,\,\,\,X(\sigma^{n-1})\in {\cal O}^n)\geq\eta_0.\label{eq:pos-hit}\end{equation}
\end{lemma}

\begin{proof}
The convergence in $(\ref{eq:converg})$ follows from time 
change and 
compactness arguments similar, for instance, to those used 
in Lemma 4.5 and Theorem 4.1 of \cite{CK18} and from the fact that $
\bar {X}$ is 
the unique natural solution of the constrained martingale problem for 
$(\frac 12\Delta ,{\cal K},B,\bar{\Xi })$ (Theorem \ref{th:cone-exist}). 

In order to prove $(\ref{eq:pos-hit})$, fix 
$\bar {x}^0\in {\cal O}\cap {\cal K}\cap\partial B_{2\delta^{*}}(
0)$, $\epsilon <\delta^{*}$ such that $B_{\epsilon}(\bar {x}^0)\subseteq 
{\cal O}\cap {\cal K}$ and $t_0>0$, 
and let $\zeta :[0,t_0]\rightarrow\R^d$ be a continuous function such that 
\begin{eqnarray*}
\zeta (0)=\bar {x},\qquad\zeta (\frac {t_0}2)=\frac {\bar {x}^0}2,\qquad 
|\zeta (s)|=\delta^{*},\;\frac {\zeta (s)}{|\zeta (s)|}\in {\cal S},\mbox{\rm \ for }
0<s<\frac {t_0}2,\\
\frac {\zeta (s)}{|\zeta (s)|}=\frac {\bar {x}^0}{|\bar {x}^0|},\quad 
|\zeta (s)|=\delta^{*}\frac 2{t_0}(t_0-s)+(2\delta^{*}+\frac {\epsilon}
2)\frac 2{t_0}(s-\frac {t_0}2),\mbox{\rm \ for }\frac {t_0}2<s\leq 
t_0.\end{eqnarray*}
Then 
\begin{eqnarray*}
&&\liminf_{n\rightarrow\infty}\P\big(\sigma^{n-1}<\vartheta,\,\,\,
X^{x^n}(\sigma^{n-1})\in {\cal O}^n\big)\\
&&\geq\liminf_{n\rightarrow\infty}\P\big(\sup_{0\leq s\leq t_0}\big
|2^{2n}\,X^{x^n}(2^{-4n}s)-\zeta (s)\big|<\frac {\epsilon}2\big)\\
&&\geq\P\big(\sup_{0\leq s\leq t_0}\big|\bar {X}^{\bar {x}}(s)-\zeta (
s)\big|<\frac {\epsilon}2\big)\\
&&=\P\big(\sup_{0\leq s\leq t_0}\big|\xi^{\bar {x}}(s)-\zeta (
s)\big|<\frac {\epsilon}2\big),\end{eqnarray*}
where $\xi$ is the absorbed process in the proof of Theorem \ref{th:cone-exist} and has a uniquely determined distribution. 
Then the assertion follows from Theorem 3.1 of 
\cite{KW91} and the Feller property of $\xi$.
\end{proof}

\begin{lemma}\label{th:coupling}
Let $\{E_n\}$ be defined by $(\ref{eq:Dn})$ and $\{Q_n\}$ be given by 
$(\ref{eq:Qk})$ and $(\ref{eq:ballcycles})$. 
For $x,\tilde {x}\in E_n$, let 
\[\epsilon_n(x,\tilde {x}):=\int\big(\tilde {f}_n(x,y)\wedge f_n(\tilde {
x},y)\big)\big(Q_n(x,dy)+Q_n(\tilde {x},dy)\big),\quad n\geq 1,\]
where 
\[\tilde {f}_n(x,\cdot ):=\frac {dQ_n(x,\cdot )}{d(Q_n(x,\cdot )+Q_
n(\tilde {x},\cdot ))},\quad f_n(\tilde {x},\cdot ):=\frac {dQ_n(\tilde {
x},\cdot )}{d(Q_n(x,\cdot )+Q_n(\tilde {x},\cdot ))}.\]
Then there exists $\epsilon_0>0$ such that 
\[\inf_{n\geq 1}\inf_{x,\tilde {x}\in E_n}\epsilon_n(x,\tilde {x}
)\geq\epsilon_0.\]
\end{lemma}

\begin{proof}
By Lemma \ref{th:scaling} and Lemma 5.3 
in \cite{CK18}, we can construct, on the same probability 
space, 
two strong Markov, natural solutions, $X$ and $\tilde {X}$, of the constrained 
martingale problem for $(A,D,B,\Xi )$, starting at $x^{}$ and $\tilde {
x}$ respectively, 
such that, denoting by $\tilde{\tau}^{n-1}$ and $\tilde{\vartheta}_
0$ the analogs of $\tau^{n-1}$ 
and $\vartheta$ for $\tilde {X}$, and by ${\cal E}$ the event 
\begin{eqnarray*}
&&{\cal E}:=\\
&&\{\exists t,\tilde {t},\,0\leq t<\tau^{n-1}\hskip-.03in\wedge\vartheta,0\leq\,
\tilde {t}<\tilde{\tau}^{n-1}\hskip-.03in\wedge\tilde{\vartheta}_0:X(t+s)
\hskip-.02in=\hskip-.03in\,\tilde {X}(\tilde {t}+s),\,0\leq s\leq (\tau^{n-1}\hskip-.03in\wedge\vartheta)
\hskip-.02in-\hskip-.02int\},\end{eqnarray*}
\[\P (\{\tau^{n-1}<\vartheta\}\cap {\cal E})=\P (\{\tilde{\tau}^{
n-1}<\tilde{\vartheta}_0\}\cap {\cal E})\geq\epsilon_0,\]
for some positive constant $\epsilon_0$ independent of $x$, $\tilde {
x}$ and $n$. 
This implies, for every $C\in {\cal B}(E_{n-1)})$,  
\begin{eqnarray*}
Q_n(x,C)&\leq&\P (\{\tau^{n-1}<\vartheta\}\cap \{X(\tau^{n-1})\in 
C\}\cap {\cal E})+\P (\tau^{n-1}<\vartheta)-\epsilon_0\\
&\leq&Q_n(\tilde {x},C)+Q_n(x,E_{n-1})-\epsilon_0,\end{eqnarray*}
and the same exchanging $x$ and $\tilde {x}$, so that 
\[\|Q_n(x,\cdot )-Q_n(\tilde {x},\cdot )\|_{TV}\leq Q_n(x,E_{n-1}
)\vee Q_n(\tilde {x},E_{n-1})-\epsilon_0.\]
On the other hand, we have 
\begin{eqnarray*}
&&Q_n(x,E_{n-1})\\
&&=\int_{\{y:\,\tilde {f}_n(x,y)>f_n(\tilde {x},y)\}}\big(\tilde {
f}_n(x,y)-\tilde {f}_n(x,y)\wedge f_n(\tilde {x},y)\big)\big(Q_n(
x,dy)+Q_n(\tilde {x},dy)\big)\\
&&\qquad\qquad\qquad\qquad\qquad\qquad +\int\big(\tilde {f}_n(x,y)\wedge 
f_n(\tilde {x},y)\big)\big(Q_n(x,dy)+Q_n(\tilde {x},dy)\big)\\
&&=\int_{\{y:\,\tilde {f}_n(x,y)>f_n(\tilde {x},y)\}}\big(\tilde {
f}_n(x,y)\vee f_n(\tilde {x},y)-\tilde {f}_n(x,y)\wedge f_n(\tilde {
x},y)\big)\big(Q_n(x,dy)+Q_n(\tilde {x},dy)\big)\\
&&\qquad\qquad\qquad\qquad\qquad\qquad +\epsilon_n(x,\tilde {
x})\\
&&\leq\|Q_n(x,\cdot )-Q_n(\tilde {x},\cdot )\|_{TV}+\epsilon_n(x,
\tilde {x}),\end{eqnarray*}
and the same exchanging $x$ and $\tilde {x}$. 
\end{proof}

%%%%%%%%%%%%%%%%%%%%%%%%%%%%%%%%%%%%%%%%%%%%%%
%% Single Appendix:                         %%
%%%%%%%%%%%%%%%%%%%%%%%%%%%%%%%%%%%%%%%%%%%%%%
%\begin{appendix}
%\section*{???}%% if no title is needed, leave empty \section*{}.
%\end{appendix}
%%%%%%%%%%%%%%%%%%%%%%%%%%%%%%%%%%%%%%%%%%%%%%
%% Multiple Appendixes:                     %%
%%%%%%%%%%%%%%%%%%%%%%%%%%%%%%%%%%%%%%%%%%%%%%
\begin{appendix}

\section{Results on the cone}\label{appcone}

\renewcommand{\theequation}{A.\arabic{equation}} 
\setcounter{equation}{0}

Let ${\cal K}$ be the cone in Condition \ref{D} and $\bar {g}$ be the 
vector field in Condition \ref{G}.  The reflecting 
Brownian motion in $\overline {{\cal K}}$ with direction of reflection $
\bar {g}$ has 
been studied by \cite{VW85} and \cite{Wil85}, for $d=2$, 
and by \cite{KW91} for $d\geq 3$, without assuming 
Conditions \ref{G} (iii) and (iv).  We summarize below 
the main results of \cite{VW85} and \cite{KW91}.  

If Conditions \ref{G} (iii) and (iv) are satisfied, a 
modification of Theorem \ref{th:exist} (Theorem 
\ref{th:cone-exist} below) yields that the reflecting 
Brownian motion in $\overline {{\cal K}}$ with direction of reflection $
\bar {g}$ is a 
semimartingale.  In dimension $d=2$, \cite{Wil85} proves 
that this is equivalent to the fact that the parameter 
$\alpha^{*}$ (defined in $(\ref{eq:alphad2})$ below) is strictly less 
than $1$.  In dimension $d\geq 3$, the issue of when the 
reflecting Brownian motion is a semimartingale is not 
discussed in \cite{KW91}.  We prove here one of the two 
implications, namely that the reflecting 
Brownian motion being a semimartingale implies that the 
parameter $\alpha^{*}$ (defined by \cite{KW91} as in Theorem 
\ref{th:funcd3}) is strictly less than $1$ (Theorem 
\ref{th:alpha}).  Beyond the intrinsic interest, this 
allows us to use some results for the cone ${\cal K}$ to obtain 
some of the estimates we need to prove uniqueness.   
(See the proof of Lemma \ref{th:auxfunc}.)  

Let ${\cal K}$ be a cone as in Condition \ref{D}, $\bar {g}$ be 
a vector field as in Conditions \ref{G} (i) and (ii) and 
$\Delta$ denote the Laplacian operator. In 
\cite{VW85}, \cite{Wil85} and \cite{KW91}, the reflecting 
Brownian motion in $\overline {{\cal K}}$ with direction of reflection $
\bar {g}$ is 
viewed as a solution to the following submartingale 
problem.  

\begin{definition}
A stochastic process $X$ with paths in ${\cal C}_{\overline {{\cal K}}}
[0,\infty )$ is a 
solution of the {\em submartingale problem\/} for $(\frac 12\Delta 
,{\cal K},\bar {g}\cdot\nabla ,\partial {\cal K})$, 
if there exists a filtration $\{{\cal F}_t\}$, on the space on which $
X$ 
is defined, such that $X$ is $\{{\cal F}_t\}$-adapted and 
\[f(X(t))-f(X(0))-\frac 12\int_0^t\Delta f(X(s))\,ds\]
is an $\{{\cal F}_t\}$-submartingale for all $f\in {\cal C}^2_b(\overline {
{\cal K}})$ such that $f$ is 
constant in a neighborhood of the origin and 
\[\bar {g}\cdot\nabla f\geq 0\quad\mbox{\rm on }\partial {\cal K}
-\{0\}.\]

The solution to the submartingale problem for 
$(\frac 12\Delta ,\overline {{\cal K}},\bar {g}\cdot\nabla ,\partial 
{\cal K})$ is unique if any two solutions have the 
same distribution.  

A solution $X$ is said to spend zero time at the origin if 
\[\E\big[\int_0^{\infty}{\bf 1}_{\{0\}}(X(s))ds\big]=0.\]
\end{definition}

For $d=2$, in polar coordinates, 
\begin{equation}{\cal K}=\{(r,z):\,r>0,\,0<z<\zeta \},\quad 0<\zeta 
<2\pi .\label{eq:C}\end{equation}
Let $ $$\partial_1{\cal K}:=\{(r,z):\,r>0,\,z=0\}$, $\partial_2{\cal K}
:=\{(r,z):\,r>0,\,z=\zeta \}$ and 
denote by $n^1$ and $n^2$ the unit inward normal vectors on 
$\partial_1{\cal K}$ and $\partial_2{\cal K}$.  Conditions \ref{G} (i) and (ii) reduce simply 
to 
\begin{equation}\bar {g}(x):=\left\{\begin{array}{cc}
\bar {g}^1,&\mbox{\rm for }x\in\partial_1{\cal K},\\
\bar {g}^2,&\mbox{\rm for }x\in\partial_2{\cal K}.\end{array}
\right.\label{eq:gd2}\end{equation}

\begin{theorem}\label{th:funcd2} (\cite{VW85})

\noindent Let $d=2$, and let ${\cal K}$ and $\bar {g}$ be as in $
(\ref{eq:C})$ 
and $(\ref{eq:gd2})$.  Let $\zeta_1$ and $\zeta_2$ denote the angles 
between $\bar {g}^1$ and $n^1$, and between $\bar {g}^2$ and $n^2$, respectively, 
taken to be positive if $\bar {g}^1$ $(\bar {g}^2$) points towards the origin.  
Set 
\begin{equation}\alpha^{*}:=\frac {\zeta_1+\zeta_2}{\zeta},\label{eq:alphad2}\end{equation}
and 
\begin{equation}\begin{array}{cc}
\psi_{\alpha^{*}}(z):=\cos(\alpha^{*}z-\zeta_1),&\mbox{\rm if }\alpha^{
*}\neq 0,\\
\psi^0(z):=\,-z\mbox{\rm tg}\,\,\zeta_1,&\mbox{\rm if }\alpha^{*}
=0.\end{array}
\label{eq:funcd2}\end{equation}
Then the function 
\begin{equation}\Psi (r,z):=\left\{\begin{array}{cc}
r^{\alpha^{*}}\psi_{\alpha^{*}}(z),&\mbox{\rm if }\alpha^{*}\neq 
0,\\
-\mbox{\rm ln}\,r\,+\,\psi^0(z),&\mbox{\rm if }\alpha^{*}=0,\end{array}
\right.\label{eq:func}\end{equation}
satisfies 
\begin{equation}\begin{array}{c}
\Delta\Psi =0\qquad\mbox{\rm in }{\cal K},\\
\bar {g}\cdot\nabla\Psi =0\qquad\mbox{\rm on }\partial {\cal K}-\{
0\}.\end{array}
\label{eq:eq}\end{equation}
\end{theorem}

\begin{theorem}\label{th:d2} (\cite{VW85}, \cite{Wil85})

\noindent Let $d=2$, and let ${\cal K}$, $\bar {g}^1,$ $\bar {g}^
2$, $\alpha^{*}$ be as in 
Theorem \ref{th:funcd2}.  For $\alpha^{*}<2$, for each $x\in\overline {
{\cal K}}$, 
there exists one and only one solution to the 
submartingale problem for $(\frac 12\Delta ,{\cal K},\bar {g}\cdot
\nabla ,\partial {\cal K})$ starting at $x$ 
that spends zero time at the origin and it is a strong 
Markov process and a Feller process. 
This solution is a semimartingale if 
and only if $\alpha^{*}<1$.  For $\alpha^{*}\geq 2$, for each $x\in\overline {
{\cal K}}$, there 
exists one and only one solution to the submartingale 
problem for $(\frac 12\Delta ,{\cal K},\bar {g}\cdot\nabla ,\partial 
{\cal K})$ starting at $x$, and it is 
absorbed at the origin after the first time it hits it.  
\end{theorem}

Now let $d\geq 3$.  Recall that $v^r$ denotes the unit radial 
vector, i.e.  $v^r(x):=\frac x{|x|}$, and $n^{{\cal K}}(x)$ denotes the unit 
inward normal to $\overline {{\cal K}}$ at $x\neq 0$.  

\begin{theorem}\label{th:funcd3} (\cite{KW91})

\noindent Let $d\geq 3$, ${\cal K}$ be a cone as in Condition \ref{D}, 
$\bar {g}$ be a vector field as in Conditions \ref{G} 8i) and (ii). 
For each $\alpha\in\R$, there exist $\lambda (\alpha )\in\R$ and 
$\psi_{\alpha}\in {\cal C}^2(\overline {{\cal S}})$ such that 
\begin{equation}\begin{array}{c}
\lambda (\alpha )\,\psi_{\alpha}+\Delta_{S^{d-1}}\psi_{\alpha}=0\qquad\mbox{\rm in }
{\cal S},\\
\alpha\bar {g}_r\psi_{\alpha}+\bar {g}_T\cdot\nabla_{S^{d-1}}\psi_{
\alpha}=0\qquad\mbox{\rm on }\partial {\cal S},\end{array}
\label{eq:funcd3}\end{equation}
where $\bar {g}_rv^r$ and $\bar {g}_T$ are the radial component and the 
component tangential to $S^{d-1}$ of $\bar {g}$.  $\psi_{\alpha}$ is strictly 
positive.  $\lambda$ and $\alpha\mapsto\psi_{\alpha}\in {\cal C}^
2(\overline {{\cal S}})$ are analytic functions.  $\lambda$ 
is concave, $\lambda (0)=0$ and 
\[\lambda'(0)=-\int_{\partial {\cal S}}\frac 1{\bar {g}\cdot n^{{\cal K}}}
\bar {g}_r\psi^{*},\]
where $\psi^{*}$ is the unique solution of 
\begin{equation}\begin{array}{c}
\Delta_{S^{d-1}}\psi^{*}=0\qquad\mbox{\rm in }{\cal S},\\
n^{{\cal K}}\cdot\nabla_{S^{d-1}}\psi^{*}-\mbox{\rm div}_{\partial 
{\cal S}}\big(\psi^{*}(\frac 1{\bar {g}\cdot n^{{\cal K}}}\bar {g}_
T-n^{{\cal K}})\big)=0\qquad\mbox{\rm on }\partial {\cal S},\end{array}
\label{eq:adjeq}\end{equation}
such that $\psi^{*}$ is strictly positive and $\int_{{\cal S}}\psi^{
*}=1$.  

If $\lambda'(0)\neq d-2$, there exists a unique $\alpha^{*}\neq 0$ such that 
\[\lambda (\alpha^{*})=\alpha^{*}(\alpha^{*}+d-2),\]
and the function $\Psi$ defined as in $(\ref{eq:func})$ for 
$\alpha^{*}\neq 0$ satisfies $(\ref{eq:eq})$.  $\alpha^{*}>0$ if $
\lambda'(0)>d-2$, 
$\alpha^{*}<0$ if $\lambda'(0)<d-2$.  

If $\lambda'(0)=d-2$, there exists a solution $\psi^0\in {\cal C}^
2(\overline {{\cal S}})$ to 
\[\begin{array}{c}
-(d-2)+\Delta_{S^{d-1}}\psi^0=0,\qquad\mbox{\rm in }{\cal S},\\
-\bar {g}_r+\bar {g}_T\cdot\nabla_{S^{d-1}}\psi^0=0,\qquad\mbox{\rm on }
\partial {\cal S},\end{array}
\]
and the function $\Psi$ defined as in $(\ref{eq:func})$ for 
$\alpha^{*}=0$ satisfies $(\ref{eq:eq})$.  In this case, we set 
$\alpha^{*}:=0$.  
\end{theorem}

\begin{theorem}\label{th:d3} (\cite{KW91})

\noindent Let $d\geq 3$ and let ${\cal K}$ and $\bar {g}$ be as in Theorem 
\ref{th:funcd3}.  For $\alpha^{*}<2$, for each $x\in\overline {{\cal K}}$, there exists 
a unique solution to the submartingale problem for 
$(\frac 12\Delta ,{\cal K},\bar {g}\cdot\nabla ,\partial {\cal K}
)$, starting at $x$, that spends zero time at 
the origin and it is a strong Markov process and a 
Feller process. 
For $\alpha^{*}\geq 2$, for each $x\in\overline {{\cal K}}$, there exists a unique solution to 
the submartingale problem for $(\frac 12\Delta ,{\cal K},\bar {g}
\cdot\nabla ,\partial {\cal K})$ starting at 
$x$, and it is absorbed at the origin after the first time it 
hits it.  
\end{theorem}

\begin{remark}
In the case $\alpha^{*}=0$, both for $d=2$ and $d\geq 3$, the function 
used in \cite{KW91} is actually $-\Psi$, but we prefer to 
have $\Psi (r,z)\rightarrow_{r\rightarrow 0}\infty$.  
\end{remark}

\begin{remark}\label{re:Holder} 
For $d=2$, $\psi_{\alpha^{*}},\psi^0\in {\cal C}^{\infty}(\overline {
{\cal S}})$.  For $d\geq 3$, a careful 
inspection of the proofs of \cite{KW91} shows that 
$\psi_{\alpha^{*}},\psi^0\in {\cal C}^{2+\beta}(\overline {{\cal S}}
)$ for every $0<\beta <1$ (see Theorem 6.31 of 
\cite{GT83} and the Remark following it.)  
\end{remark}

The function defined in $(2.7)$ of \cite{KW91}:
\begin{equation}\Phi (x):=\left\{\begin{array}{ll}
\Psi (x)^{-1},&\mbox{\rm if }\alpha^{*}<0,\\
\mbox{\rm $e^{-\Psi (x)}$},&\mbox{\rm if }\alpha^{*}=0,\\
\Psi (x),&\mbox{\rm if }\alpha^{*}>0,\end{array}
\right.\label{eq:Phi}\end{equation}
gives a way of measuring the distance from the origin 
and satisfies 
\[\bar {g}\cdot\nabla\Phi =0\mbox{\rm \ on }\partial {\cal K}-\{0
\}.\]
$\Psi$ and $\Phi$ will be used both to localize and to construct auxiliary functions 
(see Appendix \ref{appaux}). 

Let $\bar {G}$, $\bar{\Xi}$ and $B$ be as in $(\ref{eq:cone-cmp})$. 
${\cal K}$ is unbounded, but the definitions of constrained 
martingale problem, controlled martingale problem and 
natural solution of the constrained martingale problem 
carry over to $(\frac 12\Delta ,{\cal K},B,\bar{\Xi })$ without any modification, 
except that ${\cal C}^2(\overline D)$ is replaced by ${\cal C}^2_
b(\overline {{\cal K}})$ everywhere. 

\begin{lemma}\label{th:cone-bdryfunc} 
There exists a function $F\in {\cal C}^2_b(\overline {{\cal K}})$ such that 
\[\inf_{x\in\partial {\cal K},\,g\in\bar {G}(x),\,|g|=1}\nabla F(
x)\cdot g:=c_F>0.\]
\end{lemma}

\begin{proof}
See Appendix \ref{appaux}.  
\end{proof}

\begin{theorem}\label{th:cone-exist} 
Let ${\cal K}$ and $\bar {g}$ be as in Theorem \ref{th:funcd2}, for $
d=2$, 
and as in Theorem \ref{th:funcd3}, for $d\geq 3$, and, in 
addition, assume that Conditions \ref{G} (iii) 
and (iv) are satisfied. 

Then, for each $\nu\in {\cal P}(\overline {{\cal K}})$, there exists one and only one natural 
solution, $X$, to the constrained martingale problem for 
$(\frac 12\Delta ,{\cal K},B,\bar{\Xi })$ with initial distribution $
\nu$ and it is the unique 
solution of the submartingale problem for $(\frac 12\Delta ,{\cal K}
,\bar {g}\cdot\nabla ,\partial {\cal K})$ 
that spends zero time at the origin. $X$ is a 
semimartingale and a strong Markov process. The associated random measure, $\Lambda$, 
satisfies $\E[\Lambda ([0,t]\times\bar{\Xi })]<\infty$ for all $t
\geq 0$, and, for every $f\in {\cal C}^2_
b(\overline {{\cal K}})$, 
$(\ref{eq:cmp})$ is a martingale. 
\end{theorem}

\begin{proof}
Let $\{\delta_k\}$ be a strictly decreasing sequence of positive 
numbers converging to zero, with $\delta_1<1$, and let $\{D^k\}$ be 
a sequence of domains with ${\cal C}^1$ boundary such that 
$D^k\subset D^{k+1}\subset {\cal K}$, $\overline {D^k}\cap B_{\delta_
k}(0)^c=\overline {{\cal K}}\cap B_{\delta_k}(0)^c$ and 
$\overline {D^k}\cap B_{\delta_k}(0)\subset D^{k+1}.$ Also let $g^
k:\R^d-\{0\}\rightarrow\R^d$, $k\in\N$, be a 
locally Lipschitz vector field, of unit length on $\partial D^k$, such 
that $g^k(x)=\bar {g}(x)$ for $x\in\partial {\cal K}\cap B_{\delta_
k}(0)^c$ and that, denoting 
by $n^k(x)$ the unit, inward normal at $x\in\partial D^k$, it holds 
$\inf_{x\in\partial D^k}g^k(x)\cdot n^k(x)>0$.  

For each $k$, consider a sequence of bounded domains 
$\{D^{k,N}\}$, $N\in\N$, with ${\cal C}^1$ boundary, such that 
$D^{k,N}\subset D^{k,N+1}\subset D^k$, $\overline {D^{k,N}}\cap\overline {
B_N(0)}=\overline {D^k}\cap\overline {B_N(0)}$ and 
$\overline {D^{k,N}}\cap\bigg(\overline {B_N(0)}\bigg)^c\subset D^{
k,N+1}.$ Also let $g^{k,N}$ be a locally 
Lipschitz vector field, of unit length on $\partial D^{k,N}$, such that 
$g^{k,N}(x)=g^k(x)$ for $x\in\partial D^k\cap\overline {B_N(0)}$ and that, denoting by 
$n^{k,N}(x)$ the unit, inward normal at $x\in\partial D^{k,N}$, it holds 
$\inf_{x\in\partial D^k}g^{k,N}(x)\cdot n^{k,N}(x)>0$.  

Let $\xi_0$ be a random variable with compact support 
$\mbox{\rm supp}(\xi_0)\subset\overline {{\cal K}}-\{0\}$. For $k$ and $
N$ large enough that 
${\rm s}{\rm u}{\rm p}{\rm p}(\xi_0)\subset\overline {D^{k,N}}$, let $
\xi^{k,N}$ 
be the strong solution of $(\ref{eq:SDER})$ in $\overline {D^{k,N}}$ with 
direction of reflection $g^{k,N}$ and initial condition $\xi_0$, and 
let $l^{k,N}$ be the corresponding nondecreasing process.  
Define 
\[\Theta^{k,N}:=\inf\{t\geq 0:\,|\xi^{k,N}(t)|\geq N\}.\]
Let $\varphi\in {\cal C}^2(\overline {{\cal K}})$ be defined by:  
\begin{equation}\varphi (x):=\chi (\Phi (x)),\label{eq:ccphi}\end{equation}
where $\Phi$ is defined in $(\ref{eq:Phi})$ and 
$\chi :\R_{+}\rightarrow\R_{+}$ is a smooth, nondecreasing function 
such that $\chi (u)=0$, for $u\leq\sup_{|x|\leq 1}\Phi (x)$ and $
\chi (u)=u$, for 
$u\geq\inf_{|x|\geq\delta}\Phi (x)$, for $\delta>1$ such that 
$0<\sup_{|x|\leq 1}\Phi (x)<\inf_{|x|\geq\delta}\Phi (x)$.  Then 
\[\nabla\varphi (x)\cdot g^{k,N}(x)=0,\mbox{\rm \ for }x\in\partial 
D^{k,N},\,|x|\leq N,\]
\[\lim_{|x|\rightarrow\infty}\varphi (x)=\infty ,\quad\Delta\varphi 
(x)\leq c\big(1+\varphi (x)\big),\mbox{\rm \ for }x\in {\cal K},\,\,\,
|x|\leq N,\]
for some $c>0$, and, by applying Ito's formula to $\varphi$, we obtain 
\begin{equation}\lim_{N\rightarrow\infty}\sup_{k:\,{\rm supp}(\xi_0)\subset\overline {
D^k}}\P\big(\Theta^{k,N}\leq t\big)=0.\label{eq:uNnexp}\end{equation}
From $(\ref{eq:uNnexp})$, by a standard argument, we see 
that, for $\mbox{\rm supp}(\xi_0)\subset\overline {D_k}$, there is one and only one 
strong solution, $\xi^k$, to $(\ref{eq:SDER})$ in $\overline {D^k}$ with direction of 
reflection $g^k$ and initial condition $\xi_0$, and it is defined 
for all times.  Moreover, setting 
\[\theta^k:=\inf\{t\geq 0:\,\xi^k(t)\in\partial D^k\cap\overline {
B_{\delta_k}(0)}\},\quad\theta :=\lim_{k\rightarrow\infty}\theta^
k,\]
$(\ref{eq:uNnexp})$ yields that 
\begin{equation}\P\big(\theta <\infty ,\sup_{k:\,\mbox{\rm supp}(
\xi_0)\subset\overline {D^k}}\sup_{t\leq\theta^k}|\xi^k(t)|=\infty\big
)=0.\label{eq:unknexp}\end{equation}
Hence, as in Lemma \ref{th:absorb}, we can define a pair 
of stochastic processes $\xi$ and $l$ that satisfies 
$(\ref{eq:abs})$ for $0\leq t<\theta$, and almost every path of $
\xi$ 
such that $\theta <\infty$ is bounded.  For each $N\in\N$, for each 
path such that $\theta <\infty$ and $\sup_{t<\theta}|\xi (t)|\leq 
N$, we can repeat 
the argument of Lemma \ref{th:absorb} and obtain 
$(\ref{eq:xiuc})$, so that $(\ref{eq:xiuc})$ holds for almost 
every path such that $\theta <\infty$.  Therefore the solution of 
$(\ref{eq:abs})$ is well defined, up to $\theta$ included if $\theta$ is 
finite.  

We can now proceed as in Theorem 
\ref{th:exist} and construct a sequence $\{(Y^n,\lambda_0^n,\Lambda_
1^n)\}$ such 
that, for each $f\in {\cal C}^2_b(\overline {{\cal K}})$, 
\[f(Y^n(t))-f(Y^n(0))-\frac 12\int_0^t\Delta f(Y^n(s))d\lambda_0^
n(s)-\int_{[0,t]\times U}B^nf(Y^n(s),u)\Lambda^n_1(ds\times du)\]
is a martingale with respect to $\{{\cal F}^{Y^n,\lambda_0^n,\Lambda_
1^n}_t\}$, where 
\[B^nf(x,u):=u\cdot\nabla f(x){\bf 1}_{\partial {\cal K}-\{0\}}(x
)+(\rho_n)^{-1}[f(x+\rho_nu)-f(x)]{\bf 1}_{\{0\}}(x),\]
and the law of $Y^n(0)$ has compact 
support in $\overline {{\cal K}}-\{0\}$ and converges weakly to $
\nu$.  By employing again the function $\varphi$ 
defined in $(\ref{eq:ccphi})$, we can see that $\{Y^n\}$ satisfies the compact 
containment condition.  Then the same relative 
compactness arguments as in Theorem \ref{th:exist} 
apply and any limit point of $\{(Y^n,\lambda_0^n,\Lambda_1^n)\}$ is a 
solution of the controlled martingale problem for 
$(\frac 12\Delta ,{\cal K},B,\bar{\Xi })$ with initial distribution $
\nu$.  

Lemma 3.1 in \cite{CK19} holds for 
non compact state spaces as well, provided $f$ and 
$Af$ in its statement are bounded, and Lemma 
\ref{th:cone-bdryfunc} ensures that its assumptions 
are verified. In addition Lemma \ref{th:takeoff} carries over to the present context. Therefore, by 
Corollary 4.12 a) of \cite{CK19}, for each $\nu\in {\cal P}(\overline {{\cal K}})$ 
there exists a strong Markov, natural solution of the 
constrained martingale problem for $(\frac 12\Delta ,{\cal K},B,\bar{\Xi })$ with initial distribution $\nu$. 
Moreover, it can be easily checked, in the same way as in Remark 
\ref{re:clmpcont}, that all solutions to the controlled 
martingale problem for $(\frac 12\Delta ,{\cal K},B,\bar{\Xi })$ are continuous and 
that Remark \ref{re:cmpcont} and Lemma \ref{th:occtime} 
carry over to the present 
context. Thus every natural solution of the 
constrained martingale problem for $(\frac 12\Delta ,{\cal K},B,\bar{
\Xi })$ 
is a solution of the submartingale 
problem for $(\frac 12\Delta ,{\cal K},\bar {g}\cdot\nabla ,\partial 
{\cal K})$ that spends zero time at the 
origin, hence there is only one solution, $X$. 

Finally, by Lemma 3.1 of \cite{CK19}, the random measure associated to $X$, $\Lambda$, 
satisfies $\E[\Lambda ([0,t]\times\bar{\Xi })]<\infty$ for all $t
\geq 0$, and, for every $f\in {\cal C}^2_
b(\overline {{\cal K}})$, 
$(\ref{eq:cmp})$ is a martingale. 
Since the function 
$f(x):=x_i$, $i=1,...,d$, can be approximated, uniformly over 
compact sets, by functions in ${\cal C}^2_b(\overline {{\cal K}})$, $X$
is a semimartingale. 

\end{proof}

\begin{theorem}\label{th:alpha} 
Let ${\cal K}$ and $\bar {g}$ be as in Theorem \ref{th:funcd2}, for $
d=2$, 
and as in Theorem \ref{th:funcd3}, for $d\geq 3$, and, in 
addition, assume that Conditions \ref{G} (iii) 
and (iv) are satisfied. Let $\alpha^{*}$ be the parameter defined in Theorems 
\ref{th:funcd2} and \ref{th:funcd3}. Then it holds 
\[\alpha^{*}<1.\]
\end{theorem}

\begin{proof}\ 
In dimension $d=2$, the assertion follows immediately 
from Theorems \ref{th:cone-exist} and \ref{th:d2}.  By 
adapting, in a nontrivial way, an argument of 
\cite{Wil85}, we are able to prove that it holds in 
dimension $d\geq 3$ as well.  

By Theorems \ref{th:cone-exist} and \ref{th:d3}, $\alpha^{*}<2$.  
Suppose, by contradiction, that $1\leq\alpha^{*}<2$.  

Le $n^{{\cal K}}$ be the inward, unit normal to $\overline {{\cal K}}$, $
v^r(x):=\frac x{|x|}$. 
In the following, it is convenient to normalize $\bar {g}$ so that 
\[\bar {g}(x)\cdot n^{{\cal K}}(x)=1,\]
$ $rather than $|\bar {g}(x)|=1$.  Of course this does 
not affect equation $(\ref{eq:funcd3})$ and Condition 
\ref{G}.  It can be easily checked that, for $\epsilon >0$ less 
than a threshold determined by the data of the problem, 
the vector 
\begin{equation}g^{\epsilon}(z):=\bar {g}(z)-\epsilon\,v^r,\quad 
z\in\partial {\cal S},\label{eq:pertg}\end{equation}
satisfies all points of Condition \ref{G}.  Then, by 
Theorems \ref{th:cone-exist} and \ref{th:d3}, $\alpha^{\epsilon *}$, defined 
as in Theorem \ref{th:funcd3} with $\bar {g}$ replaced by 
$g^{\epsilon}$, satisfies 
\begin{equation}\alpha^{\epsilon *}<2.\label{eq:pertalpha2}\end{equation}
Let us show that 
\begin{equation}\alpha^{\epsilon *}>\alpha^{*}.\label{eq:pertalpha}\end{equation}
$ $For $\alpha >0$, let $(\lambda (\alpha ),\psi_{\alpha})$ be as in Theorem \ref{th:funcd3} 
and $(\lambda^{\epsilon}(\alpha ),\psi^{\epsilon}_{\alpha})$ be the corresponding objects with $
\bar {g}$ 
replaced by $g^{\epsilon}$.  With the notation of $(\ref{eq:funcd3})$, 
since $g^{\epsilon}_T=\bar {g}_T$, and $g^{\epsilon}_r=\bar {g}_r
-\epsilon$, $(\lambda^{\epsilon}(\alpha ),\psi^{\epsilon}_{\alpha}
)$ satisfies 
\begin{equation}\begin{array}{c}
\lambda^{\epsilon}(\alpha )\,\psi^{\epsilon}_{\alpha}+\Delta_{S^{
d-1}}\psi^{\epsilon}_{\alpha}=0,\qquad\mbox{\rm in }{\cal S},\\
\alpha (\bar {g}_r-\epsilon )\psi^{\epsilon}_{\alpha}+\bar {g}_T\cdot
\nabla_{S^{d-1}}\psi^{\epsilon}_{\alpha}=0,\qquad\mbox{\rm on }\partial 
{\cal S}.\end{array}
\label{eq:pertfuncd3}\end{equation}
Consider the function $\psi_{\alpha}/\psi^{\epsilon}_{\alpha}$.  Straightforward 
computations show that $(\ref{eq:funcd3})$ and 
$(\ref{eq:pertfuncd3})$ imply that, for $z\in {\cal S}$, 
\[\Delta_{S^{d-1}}\hskip-.02in(\psi_{\alpha}/\psi^{\epsilon}_{\alpha})(z)\hskip-.02in=\hskip-.02in[\lambda^{
\epsilon}(\alpha )-\lambda (\alpha )](\psi_{\alpha}/\psi^{\epsilon}_{
\alpha})(z)-2\big[(\psi^{\epsilon}_{\alpha})^{-1}\nabla_{S^{d-1}}\hskip-.02in
(\psi_{\alpha}/\psi^{\epsilon}_{\alpha})\cdot\nabla_{S^{d-1}}\hskip-.02in(\psi^{
\epsilon}_{\alpha})\big](z),\]
and, for $z\in\partial {\cal S}$, 
\[\big(\nabla_{S^{d-1}}(\psi_{\alpha}/\psi^{\epsilon}_{\alpha})\cdot
\bar {g}_T\big)(z)=-\alpha\,\epsilon\,(\psi_{\alpha}/\psi^{\epsilon}_{
\alpha})(z).\]
Let $z^0$ be a point of global minimum for $\psi_{\alpha}(\psi^{\epsilon}_{
\alpha})^{-1}$.  If 
$z^0\in\partial {\cal S}$, since $\bar {g}_T\cdot n^{{\cal K}}(z)
=1$, it must hold 
\[\big(\nabla_{S^{d-1}}(\psi_{\alpha}/\psi^{\epsilon}_{\alpha})\cdot
\bar {g}_T\big)(z^0)\geq 0,\]
while 
\[-\alpha\,\epsilon\,(\psi_{\alpha}/\psi^{\epsilon}_{\alpha})(z^0
)<0,\qquad\forall\alpha >0,\]
because $\psi_{\alpha}$ and $\tilde{\psi}_{\alpha}$ are strictly positive.  Therefore it 
must be $z^0\in {\cal S}$ and 
\[\nabla_{S^{d-1}}(\psi_{\alpha}/\psi^{\epsilon}_{\alpha})(z^0)=0
,\quad\Delta_{S^{d-1}}(\psi_{\alpha}/\psi^{\epsilon}_{\alpha})(z^
0)>0,\]
which yields 
\begin{equation}\lambda^{\epsilon}(\alpha )>\lambda (\alpha ),\qquad
\forall\alpha >0.\label{eq:pertlambda}\end{equation}
Then $\big(\lambda^{\epsilon}\big)'(0)\geq\lambda'(0)>d-2$, so that $
\alpha{}^{\epsilon *}>0$.  Hence, taking 
into account that $\lambda (\alpha )-\alpha (\alpha +d-2)$ vanishes for $
\alpha =0$ 
and $\alpha =\alpha^{*}$ and is strictly concave, $(\ref{eq:pertlambda})$ 
gives $(\ref{eq:pertalpha})$.  

The function $\Psi^{\epsilon}$ defined by $(\ref{eq:funcd3})$ with 
$\psi_{\alpha}$ replaced by $\psi^{\epsilon}_{\alpha}$ and $\alpha^{
*}$ replaced by $\alpha^{\epsilon}{}^{*}$ has the following properties: 
\begin{eqnarray}
\label{eq:funceps}\Psi^{\epsilon\in {\cal C}^2(\overline {{\cal K}}
-\{0\})\cap {\cal C}^1(\overline {{\cal K})}},\quad\nabla\Psi^{\epsilon}
(0)=0,\qquad\qquad\non\\
\Delta\Psi^{\epsilon}(x)=0,\quad\in {\cal K},\qquad\qquad\qquad\qquad\non\\
\hfill \\
\epsilon c_1\Psi (x)^{(\alpha^{\epsilon}{}^{*}-1)/\alpha^{*}}\leq\big
(\bar {g}\cdot\nabla\Psi^{\epsilon}\big)(x)\leq\epsilon c_2\Psi (
x)^{(\alpha^{\epsilon}{}^{*}-1)/\alpha^{*}},\quad\mbox{\rm on $\partial 
{\cal K}-\{0\}$},\non\\
c_1\Psi (x)^{\alpha^{\epsilon}{}^{*}/\alpha^{*}}\leq\Psi^{\epsilon}
(x)\leq c_2\Psi (x)^{\alpha^{\epsilon}{}^{*}/\alpha^{*}},\qquad\mbox{\rm in }\overline {
{\cal K}}\non\qquad\quad\end{eqnarray}
where the constants $c_1$ and $c_2$ can be taken independent 
of $\epsilon$ because the map $\alpha\rightarrow\psi_{\alpha}$ is continuous and 
$1\leq\alpha^{\epsilon *}\leq 2$. 
Of course $(\ref{eq:funceps})$ still holds if we revert 
to the usual normalization of $\bar {g}$, $|\bar {g}|=1$, as we will do for 
the rest of the proof. 

Let $X$ be the solution of the submartingale 
problem for $(\frac 12\Delta ,{\cal K},\bar {g}\cdot\nabla ,\partial 
{\cal K})$ spending zero time at the 
origin and starting at $x=0$.  
The rest of the proof uses essentially the same 
argument as Theorem 5 of \cite{Wil85}, but our proof is 
simpler because we can take advantage of the fact that  
$X$ is the solution of the constrained martingale 
problem for $(\frac 12\Delta ,{\cal K},B,\overline {\Xi})$. Fix $
0<\delta <1$, and let 
\[T^1:=\inf\{t\geq 0:\,\Psi (X(t))\geq 1\},\]
\[X^1(t):=X(t\wedge T^1).\]
Define 
\[\vartheta^1_0:=\inf\{t\geq 0:\,X^1(t)=0\}=0,\]
\[\theta_n^1:=\inf\{t\geq\vartheta_{n-1}^1:\,\Psi (X^1(t))=\delta 
\},\quad n\geq 1,\]
\[\vartheta_n^1:=\inf\{t\geq\theta_n^1:\,X^1(t)=0\},\quad n\geq 1
,\]
with the usual convention that the infimum of the 
empty set is $\infty$. By the continuity of $X$, $\theta_n^1\uparrow
\infty$ and $\vartheta_n^1\uparrow\infty$ as $n\rightarrow\infty$. We have 
\begin{eqnarray}
&&\Psi^{\epsilon}(X^1(t))\non\\
&=&\quad\sum_{n=0}^{\infty}{\bf 1}_{\{\vartheta_n^1\leq t\}}[\Psi^{
\epsilon}(X^1(t\wedge\theta_{n+1}^1))-\Psi^{\epsilon}(X^1(t\wedge
\vartheta_n^1))]\label{eq:cycles}\\
&&\quad\,\,\,+\sum_{n=1}^{\infty}{\bf 1}_{\{\theta_n^1\leq t\}}[\Psi^{
\epsilon}(X^1(t\wedge\vartheta_n^1))-\Psi^{\epsilon}(X^1(t\wedge\theta_
n^1))].\non\end{eqnarray}
As far as the first summand is concerned, we have, by 
$(\ref{eq:funceps})$, on the set $\{\vartheta_n^1\leq t\}$, 
\[|\Psi^{\epsilon}(X^1(t\wedge\theta_{n+1}^1)-\Psi^{\epsilon}(X^1
(\vartheta_n^1))|=\Psi^{\epsilon}(X^1(t\wedge\theta_{n+1}^1)\leq 
c_2\delta^{\alpha^{\epsilon}{}^{*}/\alpha^{*}}.\]
In addition, it can be easily checked that the argument 
used to prove $(52)$ in \cite{Wil85}, combined with Lemma 
2.8 of \cite{KW91}, still works, that is 
\begin{equation}\E^0\big[\sum_{n=1}^{\infty}{\bf 1}_{\{\vartheta_
n^1\leq t\}}\big]\leq\E^0\big[\sum_{n=1}^{\infty}{\bf 1}_{\{\theta_
n^1\leq t\}}\big]\leq c\,\delta^{-1}\,\frac {t+1}{1-\delta^{2/\alpha^{
*}}},\label{eq:hitfreq}\end{equation}
where $c$ depends only on $\alpha^{*}$ and $\psi_{\alpha^{*}}$. Thus, 
for each $\epsilon$, by $(\ref{eq:pertalpha})$, 
the expectation of the first summand in $(\ref{eq:cycles})$ 
vanishes as $\delta\rightarrow 0$. 
As for the second summand, by 
$(\ref{eq:hitfreq})$ and the definition of $X^1$, 
it is bounded above by an integrable random 
variable. Moreover, taking into account that $\theta_n^1\leq t$ implies 
$\theta_n^1\leq T^1$, hence $\{\theta_n^1\leq t\}\in {\cal F}_{\theta_
n^1\wedge T^1}$, we have 
\begin{eqnarray}
&&\E^0\big[\sum_{n=1}^{\infty}{\bf 1}_{\{\theta_n^1\leq t\}}[\Psi^{
\epsilon}(X^1(t\wedge\vartheta_n^1))-\Psi^{\epsilon}(X^1(\theta_n^
1))]\big]\non\\
&=&\E^0\big[\sum_{n=1}^{\infty}{\bf 1}_{\{\theta_n^1\leq t\}}[\Psi^{
\epsilon}(X(t\wedge\vartheta_n^1\wedge T^1))-\Psi^{\epsilon}(X(\theta_
n^1\wedge T^1))]\big]\non\\
&=&\E^0\bigg[\sum_{n=1}^{\infty}{\bf 1}_{\{\theta_n^1\leq t\}}\E\big
[\Psi^{\epsilon}(X(t\wedge\vartheta_n^1\wedge T^1))-\Psi^{\epsilon}
(X(\theta_n^1\wedge T^1))\big|{\cal F}_{\theta_n^1\wedge T^1}\big
]\bigg].\non\end{eqnarray}
Recall that $\bar{\Xi}$ and $B$ are defined by $(\ref{eq:cone-cmp})$. Note  
that, by $(\ref{eq:funceps})$, $B\Psi^{\epsilon}$ is a continuous function on $
\bar{\Xi}$ 
and $0\leq B\Psi^{\epsilon}(x,u)|\leq\epsilon c_2\Psi (x)$ for all $
(x,u)\in\bar{\Xi}$, and that, by 
the definition of natural solution, 
\begin{eqnarray*}
&&\int_{\big[\theta_n^1\wedge T^1,t\wedge\vartheta_n^1\wedge T^1\big
]\times\bar{\Xi}}B\Psi^{\epsilon}(x,u)\,\Lambda (ds\times dx\times 
du)\\
&&=\int_{\big[\theta_n^1\wedge T^1,t\wedge\vartheta_n^1\wedge T^1\big
]\times\bar{\Xi}}B\Psi^{\epsilon}(X(s),u)\,\Lambda (ds\times dx\times 
du)\\
&&=\int_{\big[\theta_n^1\wedge T^1,t\wedge\vartheta_n^1\wedge T^1\big
]\times\bar{\Xi}}\big(B\Psi^{\epsilon}(x,u)\wedge\epsilon c_2\big
)\,\Lambda (ds\times dx\times du)\end{eqnarray*}
Then, taking into account that $\E[\Lambda ([0,t]\times\bar{\Xi }
)]<\infty$ 
(Theorem \ref{th:cone-exist}), although $\Psi^{\epsilon}\notin {\cal C}^
2_b(\overline {{\cal K}})$ the above chain of 
equalities can be continued as 
\begin{eqnarray*}
&&=\E^0\bigg[\sum_{n=1}^{\infty}{\bf 1}_{\{\theta_n^1\leq t\}}\E\bigg
[\int_{\big[\theta_n^1\wedge T^1,t\wedge\vartheta_n^1\wedge T^1\big
]\times\bar{\Xi}}\big(B\Psi^{\epsilon}(x,u)\wedge\epsilon c_2\big
)\Lambda (ds\times dx\times du)\big|{\cal F}_{\theta_n^1\wedge T^
1}\bigg]\bigg]\non\\
&&=\E^0\bigg[\sum_{n=1}^{\infty}{\bf 1}_{\{\theta_n^1\leq t\}}\int_{\big
[\theta_n^1\wedge T^1,t\wedge\vartheta_n^1\wedge T^1\big]\times\bar{
\Xi}}\big(B\Psi^{\epsilon}(x,u)\,\wedge\epsilon c_2\big)\Lambda (
ds\times dx\times du)\bigg]\non\\
&&\leq\E^0\bigg[\int_{\big[0,t\wedge T^1\big]\times\bar{\Xi}}\big
(B\Psi^{\epsilon}(X(s),u)\wedge\epsilon c_2\big)\Lambda (ds\times 
dx\times du)\bigg].\non\end{eqnarray*}
Summing up, we have proved that, for each $\epsilon$, 
\[\E^0\big[\Psi^{\epsilon}(X(t\wedge T^1))\big]\leq\E^0\bigg[\int_{\big
[0,t\wedge T^1\big]\times\bar{\Xi}}\big(B\Psi^{\epsilon}(x,u)\wedge
\epsilon c_2\big)\Lambda (ds\times dx\times du)\bigg].\]
By $(\ref{eq:funceps})$, for $\Psi (x)\leq 1$, $\Psi^{\epsilon}(x
)\geq c_1\Psi (x)^2$, so 
that, in the limit as $\epsilon\rightarrow 0$, we find 
\[\E^0\big[\Psi (X(t\wedge T^1))^2\big]=0,\quad\forall t>0,\]
which contradicts the fact that $X$ spends zero time at 
the origin. 
\end{proof}

\section{Auxiliary functions}\label{appaux}
\renewcommand {\theequation}{B.\arabic{equation}} 
\setcounter{equation}{0} 

\begin{lemma}\label{th:loc} 
There exists $\delta^{*}>0$ and a function $V\in {\cal C}^2(\overline 
D-\{0\})$ such 
that 
\begin{equation}V(x)>0,\quad\mbox{\rm for }x\in\overline D-\{0\},
\quad\lim_{x\in\overline D,\,x\rightarrow 0}V(x)=0,\label{eq:loc}\end{equation}
\begin{equation}\nabla V(x)\cdot g(x)\leq 0,\quad\mbox{\rm for }x
\in\big(\partial D-\{0\}\big)\cap\overline {B_{\delta^{*}}(0)}.\label{eq:nablaloc}\end{equation}
Of course we can always define $V(0):=0$.
\end{lemma}

\begin{proof}
It is enough to prove that there exists $\delta^{*}>0$ and a function $
V\in {\cal C}^2\big(\big((\overline D-\{0\}\big)\cap\overline {B_{
\delta^{*}}(0)}\big)$ such 
that $V(x)>0$ for $x\in\big(\overline D-\{0\}\big)\cap\overline {
B_{\delta^{*}}(0)}$, the limit in 
$(\ref{eq:loc})$ holds and $(\ref{eq:nablaloc})$ is satisfied. 

Let $\alpha^{*}$, $\psi_{\alpha^{*}}$ and $\psi^0$ be as in Theorem \ref{th:funcd2}, 
for $d=2$, and in Theorem \ref{th:funcd3}, for $d\geq 3$. 
Since $\partial {\cal S}$ is smooth, by Condition \ref{D} (ii), we can extend 
$\psi_{\alpha^{*}}$ to a ${\cal C}^2$ function on some open neighborhood 
${\cal S}^{*}$ of $\overline {{\cal S}}$ such that 
\[\inf_{z\in {\cal S}^{*}}\psi_{\alpha^{*}}(z)>0.\]
Analogously we can extend $\psi^0$ to a ${\cal C}^2$ function on some open neighborhood 
${\cal S}^{*}$ of $\overline {{\cal S}}$ such that $\inf_{z\in {\cal S}^{
*}}e^{-\psi^0(z)}>0$. Let ${\cal K}^{*}:=\{x:\,x=rz,\,z\in {\cal S}^{
*},\,r>0\}$. 

Let $\Phi$ be the function defined in $(\ref{eq:Phi})$. We have 
\begin{equation}\Phi (x)>0\;\mbox{\rm \ in }K^{*},\quad\lim_{x\in 
{\cal K}^{*},\,x\rightarrow 0}\Phi (x)=0,\quad\label{eq:Phizero}\end{equation}
and, for some $c_{\Phi}\geq c'_{\Phi}>0$, 
\begin{equation}\frac {c_{\Phi}'\Phi (x)}{|x|}\leq |\nabla\Phi (x
)|\leq\frac {c_{\Phi}\Phi (x)}{|x|},\quad\,|D^2\Phi (x)|\leq\frac {
c_{\Phi}\Phi (x)}{|x|^2},\quad x\in {\cal K}^{*}.
\label{eq:Phiderb}\end{equation}
We will look for a function $V$ of the form 
\begin{equation}V(x):=f(\Phi (x))-c_Ve\cdot x,\label{eq:Vshape}\end{equation}
for some $f\in {\cal C}^2((0,\infty ))$ such that $\lim_{u\rightarrow 
0^{+}}f(u)=0$ and some 
$c_V>0$. Then, by Condition \ref{G}(iv), 
\[\nabla V(x)\cdot\bar {g}(x)\leq -c_Vc_e,\quad\mbox{\rm for }x\in
\partial {\cal K}-\{0\}.\]
By Condition \ref{D} (i), there is $\delta^{*}$, $0<\delta^{*}\leq 
r_D$, such that 
$\big(\overline D-\{0\}\big)\cap\overline {B_{\delta^{*}}(0)}\subset 
{\cal K}^{*}\cap\overline {B_{\delta^{*}}(0)}$. Then, for 
$x\in\big(\partial D-\{0\}\big)\cap\overline {B_{\delta^{*}}(0)}$, letting $
\bar {z}$ be the closest point on 
$\partial {\cal S}$ to $\frac x{|x|}$, by Condition \ref{D} (i) and Condition \ref{G} 
(i), we have 
\begin{eqnarray*}
&&\nabla V(x)\cdot g(x)\\
&&\leq\nabla V(|x|\bar {z})\cdot\bar {g}(|x|\bar {z})\,+\,|\nabla 
V(x)-\nabla V(|x|\bar {z})|\,+\,|\nabla V(|x|\bar {z})|\,|g(x)-\bar {
g}(|x|\bar {z})|\\
&&\leq -c_Vc_e\,+\,d\,\sup_{0<t<1}\big|D^2V(tx+(1-t)|x|\bar {z})\big
|c_D\,|x|^2\,+\,|\nabla V(|x|\bar {z})|c_g\,|x|.\end{eqnarray*}
 Since $\inf_{0<t<1}\big|tx+(1-t)|x|\bar {z}\big|\geq\frac 12|x|$ for $
|x|\leq\delta^{*}$, 
$\delta^{*}\leq\sqrt {3}/c_D$, one way to ensure $(\ref{eq:nablaloc})$, 
is to choose $f$ in $(\ref{eq:Vshape})$ so that 
\[\lim_{x\in {\cal K}^{*},\,x\rightarrow 0}\,|\nabla V(x)|\,|x|=0
,\quad\lim_{x\in {\cal K}^{*},\,x\rightarrow 0}\,|D^2V(x)|\,|x|^2
=0,\]
that is 
\[\begin{array}{c}
\lim_{x\in {\cal K}^{*},\,x\rightarrow 0}\,|f'(\Phi (x))|\,|\nabla
\Phi (x)|\,|x|=0,\\
\\
\lim_{x\in {\cal K}^{*},\,x\rightarrow 0}\,|f^{\prime\prime}(\Phi 
(x))|\,|\nabla\Phi (x)|^2\,|x|^2=0,\quad\lim_{x\in {\cal K}^{*},\,
x\rightarrow 0}|f'(\Phi (x))|\,|D^2\Phi (x)|\,|x|^2=0.\end{array}
\]
In view of $(\ref{eq:Phiderb})$, this is implied by 
\[\lim_{x\in {\cal K}^{*},\,x\rightarrow 0}f'(\Phi (x))\Phi (x)\,
=0,\lim_{x\in {\cal K}^{*},\,x\rightarrow 0}f^{\prime\prime}(\Phi 
(x))\Phi (x)^2\,=0.\]
If, in addition,   
\[\inf_{x\in\big(\overline D-\{0\}\big)\cap\overline {B_{\delta^{
*}}(0)}}\frac {f(\Phi (x))}{|x|}>0,\]
then, by choosing $c_V<\inf_{x\in\big(\overline D-\{0\}\big)\cap\overline {
B_{\delta^{*}}(0)}}\frac {f(\Phi (x))}{|x|}$, we will obtain 
$V(x)>0$ for $x\in\big(\overline D-\{0\}\big)\cap\overline {B_{\delta^{
*}}(0)}$. 

Therefore we can take, for instance,  
\[f(u):=u^{1/|\alpha^{*}|},\mbox{\rm \ for }\alpha^{*}\neq 0,\qquad 
f(u):=u,\mbox{\rm \ for }\alpha^{*}=0.\]
\end{proof}
\vskip.3in

\noindent {\bf Proof of Lemma \ref{th:bdryfunc}.}
\medskip

\noindent Let $\delta^{*}$ and $V$ be as in Lemma \ref{th:loc}. 

By Condition \ref{G} (i) and (iv), possibly 
by taking a smaller $\delta^*$, we can always suppose that 
\[\inf_{g\in G(x),|g|=1,\,x\in\partial D,\,|x|\leq\delta^{*}}e\cdot 
g>0.\]

Let $0<p^{*}<1$ be such that 
\[\sup_{x\in\overline D,\,|x|\leq p^{*}\delta^{*}}V(x)<\inf_{x\in\overline 
D,\,|x|\geq\delta^{*}}V(x).\]
Let $\tilde {D}$ be a bounded domain with ${\cal C}^1$ boundary such that $
\tilde {D}\subset D$ and 
$\overline {\tilde {D}}\cap B_{p^{*}\delta^{*}}(0)^c=\overline D\cap 
B_{p^{*}\delta^{*}}(0)^c$ and let 
$\tilde {g}:\R^d\rightarrow\R^d$ be a locally Lipschitz vector field, of 
unit length on $\partial\tilde {D}$, such that $\tilde {g}(x)=g(x
)$ for $x\in\partial\tilde {D}\cap B_{p^{*}\delta^{*}}(0)^c$ and, denoting by 
$\tilde {n}(x)$ the unit, inward normal at $x\in\partial\tilde {D}$, it holds  
$\inf_{x\in\partial\tilde {D}}\tilde {g}(x)\cdot\tilde {n}(x)>0$. 
There exists a function $\tilde {F}\in {\cal C}^2\big(\overline {\tilde {
D}}\big)$ such that 
\[\inf_{x\in\partial\tilde {D}}\nabla\tilde {F}(x)\cdot\tilde {g}
(x)>0,\]
(see \cite{CIL92}, Lemma 7.6). Of course we can 
always assume that 
\[\sup_{x\in\overline D,\,\,p^{*}\delta^{*}\leq |x|\leq\delta^{*}}
\tilde {F}(x)\leq -\delta^{*}.\]

Now let $\chi :\R\rightarrow [0,1]$ be a 
nonincreasing, ${\cal C}^{\infty}$ function such that $\chi (u)=1$ for 
$u\leq\sup_{x\in\overline D,\,|x|\leq p^{*}\delta^{*}}V(x)$ and $
\chi (u)=0$ for 
$u\geq\inf_{x\in\overline D,\,|x|\geq\delta^{*}}V(x)$ . Defining 
\[F(x):=\chi (V(x))\,e\cdot x+\big(1-\chi (V(x))\big)\tilde {F}(x
),\]
we have 
\[\nabla F(x)=\big[\chi (V(x))\,e+\big(1-\chi (V(x))\big)\nabla\tilde {
F}(x)\big]+\big(e\cdot x-\tilde {F}(x)\big)\chi'(V(x))\nabla V(x)\]
so that, for all $x\in\partial D$, $g\in G(x)$, $|g|=1$,   
\[\nabla F(x)\cdot g\geq\inf_{g\in G(y),|g|=1,\,y\in\partial D,\,
|y|\leq\delta^{*}}e\cdot g\wedge\inf_{y\in\partial\tilde {D}}\nabla
\tilde {F}(y)\cdot\tilde {g}(y).\]
\hfill $\Box$ \medskip
\vskip.3in

\noindent {\bf Proof of Lemma \ref{th:auxfunc}.}
\medskip

\noindent
Let $\alpha^{*}$, $\psi_{\alpha^{*}}$ and $\psi^0$ be as in Theorem \ref{th:funcd2}, 
for $d=2$, and in Theorem \ref{th:funcd3}, for $d\geq 3$, and let   
$\Psi$ be given by $(\ref{eq:func})$. 
By Theorem \ref{th:alpha}, we can fix $\beta$ such that  
\begin{equation}\alpha^{*}\vee 0<\beta <1.\label{eq:beta}\end{equation}
Then, by Remark \ref{re:Holder}, we can extend $\psi_{\alpha^{*}}$ to a $
{\cal C}^{2+\beta}$ function on some open neighborhood 
${\cal S}^{*}$ of $\overline {{\cal S}}$ such that 
\[\inf_{z\in {\cal S}^{*}}\psi_{\alpha^{*}}(z)>0.\]
Analogously we can extend $\psi^0$ to a ${\cal C}^{2+\beta}$ function on some open neighborhood 
${\cal S}^{*}$ of $\overline {{\cal S}}$. 
Let ${\cal K}^{*}:=\{x:\,x=rz,\,z\in {\cal S}^{*},\,r>0\}$. We will choose $
{\cal S}^{*}$ such that 
\begin{equation}e\cdot x\geq -c^{*}_e|x|,\,\,\,\,0<c^{*}_e<1,\qquad 
x\in {\cal K}^{*}.\label{eq:xpos}\end{equation}

The derivatives of $\Psi$ satisfy, for some $c_{\Psi}>c'_{\Psi}
>0$, 
\begin{equation}\frac {c_{\Psi}'}{|x|}\leq |\nabla\Psi (x)|\leq\frac {
c_{\Psi}}{|x|},\quad\,|D^2\Psi (x)|\leq\frac {c_{\Psi}}{|x|^2},\quad 
x\in {\cal K}^{*},\quad\mbox{\rm if }\alpha^{*}=0,\label{eq:Psi0derb}\end{equation}
\begin{equation}\frac {c_{\Psi}'\Psi (x)}{|x|}\leq |\nabla\Psi (x
)|\leq\frac {c_{\Psi}\Psi (x)}{|x|},\quad\,|D^2\Psi (x)|\leq\frac {
c_{\Psi}\Psi (x)}{|x|^2},\quad x\in {\cal K}^{*},\quad\mbox{\rm if }
\alpha^{*}\neq 0.\label{eq:Psiderb}\end{equation}
Let $\delta^{*}$, $0<\delta^{*}\leq r_D$, be such that $\big(\overline 
D-\{0\}\big)\cap\overline {B_{\delta^{*}}(0)}\subset {\cal K}^{*}
\cap\overline {B_{\delta^{*}}(0)}$. 
The fact that $\psi^0,\psi_{\alpha^{*}}\in {\cal C}^{2+\beta}({\cal S}^{
*})$, combined with 
$(\ref{eq:eq})$ and Condition \ref{D} (i), implies that  
\begin{equation}\,|\Delta\Psi (x)|\leq\frac {c_{\Psi}}{|x|^{2-\beta}}
,\quad x\in\big(\overline D-\{0\}\big)\cap\overline {B_{\delta^{*}}
(0)},\qquad\mbox{\rm if }\alpha^{*}=0,\label{eq:Psi0deltab}\end{equation}
\begin{equation}|\Delta\Psi (x)|\leq c_{\Psi}\frac {\Psi (x)}{|x|^{
2-\beta}},\quad x\in\big(\overline D-\{0\}\big)\cap\overline {B_{
\delta^{*}}(0)},\qquad\mbox{\rm if }\alpha^{*}\neq 0.\label{eq:Psideltab}\end{equation}
Consider first the case $\alpha^{*}\leq 0$. We look for $V$ of the 
form 
\[V(x):=f(\Psi (x))-e\cdot x,\]
with $f\in {\cal C}^2((0,\infty ))$ such that 
\[\lim_{x\in\big(\overline D-\{0\}\big),\,x\rightarrow 0}f(\Psi (
x))\,=\infty .\]
By the same computations as in the proof of Lemma 
\ref{th:loc} and by $(\ref{eq:Psi0derb})$, $(\ref{eq:Psiderb})$, 
we see that $(\ref{eq:nablaV-})$ is verified as soon as 
\begin{equation}\lim_{x\in {\cal K}^{*},\,x\rightarrow 0}f'(\Psi 
(x))\,=0,\lim_{x\in {\cal K}^{*},\,x\rightarrow 0}f^{\prime\prime}
(\Psi (x))=0,\quad\mbox{\rm if }\alpha^{*}=0,\label{eq:f0der0}\end{equation}
\begin{equation}\lim_{x\in {\cal K}^{*},\,x\rightarrow 0}f'(\Psi 
(x))\Psi (x)\,=0,\lim_{x\in {\cal K}^{*},\,x\rightarrow 0}f^{\prime
\prime}(\Psi (x))\Psi (x)^2=0,\quad\mbox{\rm if }\alpha^{*}\neq 0
.\label{eq:fder0}\end{equation}
As far as $(\ref{eq:AV-})$ is concerned, we have, 
by Condition \ref{A}, $(\ref{eq:Psi0derb})$ and 
$(\ref{eq:f0der0})$, or $(\ref{eq:Psiderb})$, and 
$(\ref{eq:fder0})$, 
\begin{eqnarray}
&&2AV(x)\non\\
&&=\,\Delta V(x)+\,|x|^{-1}o(1)\label{eq:AVcomp}\\
&&=f^{\prime\prime}(\Psi (x))\,|\nabla\Psi (x)|^2+f'(\Psi (x))\,\Delta
\Psi (x)+|x|^{-1}o(1)\non\\
&&=|x|^{\beta -2}\bigg(f^{\prime\prime}(\Psi (x))\,|\nabla\Psi (x
)|^2|x|^{2-\beta}+f'(\Psi (x))\,\Delta\Psi (x)|x|^{2-\beta}+|x|^{
1-\beta}o(1)\bigg).\non
\end{eqnarray}

Hence, taking into account $(\ref{eq:Psi0derb})$ and 
$(\ref{eq:Psi0deltab})$ or 
$(\ref{eq:Psiderb})$ and $(\ref{eq:Psideltab})$, 
$(\ref{eq:AV-})$ holds for instance if, in addition to $(\ref{eq:f0der0})$ 
or $(\ref{eq:fder0})$, 
\[\sup_{x\in (\overline D-\{0\})\cap\overline {B_{\delta^{*}}(0)}}
f^{\prime\prime}(\Psi (x))\,|x|^{-\beta}\,<0,\quad\mbox{\rm if }\alpha^{
*}=0,\]
\[\sup_{x\in (\overline D-\{0\})\cap\overline {B_{\delta^{*}}(0)}}
f^{\prime\prime}(\Psi (x))\,\Psi (x)\,^2\,|x|^{-\beta}\,<0,\quad\mbox{\rm if }
\alpha^{*}<0.\]
Therefore, for $\delta^{*}$ small enough that $\Psi (x)>1$, we can take 
\[f(u):=\ln(u),\mbox{\rm \ for }\alpha^{*}=0,\qquad f(u):=\ln(\ln
(u)),\mbox{\rm \ for }\alpha^{*}<0.\]

In the case $0<\alpha^{*}<1$, we look for $V_1$ and $V_2$ of the form 
\[V_1(x):=f_1(\Psi (x))+e\cdot x,\qquad V_2(x):=f_2(\Psi (x))-e\cdot 
x,\]
with $f_1,f_2\in {\cal C}^2((0,\infty ))$ such that 
\[\lim_{x\in\overline D,\,x\rightarrow 0}f_1(\Psi (x))\,=0,\qquad\liminf_{
x\in\overline D,\,x\rightarrow 0}f_1(\Psi (x))|x|^{-1}\,\geq 0,\]
\[\lim_{x\in\overline D,\,x\rightarrow 0}f_2(\Psi (x))\,=0,\qquad\liminf_{
x\in\overline D,\,x\rightarrow 0}f_2(\Psi (x))|x|^{-1}\,>1.\]
Consider $V_2$. If 
\begin{equation}\sup_{x\in (\overline {{\cal K}^{*}}-\{0\})\cap\overline {
B_{\delta^{*}}(0)}}|f_2'(\Psi (x))\,|\,<\infty ,\qquad\sup_{x\in 
(\overline {{\cal K}^{*}}-\{0\})\cap\overline {B_{\delta^{*}}(0)}}
|f_2^{\prime\prime}(\Psi (x))\,|\Psi (x)<\infty ,\label{eq:f+derb}\end{equation}
we have, by the same computations as for $(\ref{eq:AVcomp})$, 
\begin{eqnarray*}
&&2AV_2(x)\\
&&=\,f_2^{\prime\prime}(\Psi (x))\,|\nabla\Psi (x)|^2+f_2'(\Psi 
(x))\,\Delta\Psi (x)+|x|^{-1}\Psi (x)O(1)\\
&&=|x|^{\beta -2}\Psi (x)\bigg(f_2^{\prime\prime}(\Psi (x))\,|\nabla
\Psi (x)|^2\Psi (x)^{-1}|x|^{2-\beta}+f_2'(\Psi (x))\,\Delta\Psi 
(x)\Psi (x)^{-1}|x|^{2-\beta}\\
&&\hskip3.1in+|x|^{1-\beta}O(1)\bigg).\end{eqnarray*}
If, in addition to $(\ref{eq:f+derb})$, $f_2$ satisfies 
\[\lim_{x\in\overline D,\,x\rightarrow 0}f_2^{\prime\prime}(\Psi 
(x))\,\Psi (x)\,\,|x|^{-\beta}\,=-\infty ,\]
then, taking into account 
$(\ref{eq:Psiderb})$, $V_2$ satisfies $(\ref{eq:AV+})$. 
Note that $(\ref{eq:f+derb})$ implies 
$(\ref{eq:fder0})$, so that $(\ref{eq:nablaV+})$ is verified 
as well. By analogous computations we see that if $f_1$ 
satisfies $(\ref{eq:f+derb})$ and 
\[\lim_{x\in\overline D,\,x\rightarrow 0}f_1^{\prime\prime}(\Psi 
(x))\,\Psi (x)\,\,|x|^{-\beta}\,=\infty ,\]
then $V_1$ will satisfy $(\ref{eq:nablaV+})$ and $(\ref{eq:AV+})$. 
With the choice 
\[f_1(u)=\exp(u)-1,\qquad f_2(u)=\ln(u+1),\]
$V_1$ and $V_2$ will satisfy the third condition in $(\ref{eq:V+zero})$ as well. 
\hfill $\Box$ 
\vskip.3in

\noindent {\bf Proof of Lemma \ref{th:cone-bdryfunc}.}
\medskip
\noindent
Since $\partial {\cal S}$ is of class ${\cal C}^3$, there exists $
F^{{\cal S}}\in {\cal C}^2(\overline {{\cal S}})$ such that 
$F^{{\cal S}}(z)\geq 0$ and, on $\partial {\cal S}$, $F^{{\cal S}}(
z)=0$, $\nabla_{S^{d-1}}F^{{\cal S}}(z)\cdot\bar {g}(z)\geq 1$. 
Let $\chi^{*}\in {\cal C}^{\infty}_b([0,\infty ))$ be a nondecreasing function such that 
$\chi^{*}(u)=u$ for $u\leq 1/2$, $\chi^{*}(u)=1$ for $u\geq 1$. Define 
\[F^{*}(x):=\chi^{*}(|x|F^{{\cal S}}(\frac x{|x|})).\]
Then $F^{*}\in {\cal C}^2_b(\overline {{\cal K}}-\{0\})$ and, for $
x\in\partial {\cal K}-\{0\}$, 
\[\nabla F^{*}(x)\cdot\bar {g}(x)=\nabla_{{\cal S}^{d-1}}F^{{\cal S}}(\frac x{|x|})\cdot\bar {g}(\frac 
x{|x|})\geq c^{*\prime}.\]

Now let $\delta >1$ be such that 
$\sup_{x\in\overline {{\cal K}},\,|x|\leq 1}\Phi (x)<\inf_{x\in\overline {
{\cal K}},\,|x|\geq\delta}\Phi (x)$. 
Let $D$ be a bounded domain such that 
$D\subset {\cal K}\cap B_{\delta +1}(0)$, $\overline D\cap\overline {
B_{\delta}(0)}=\overline {{\cal K}}\cap\overline {B_{\delta}(0)}$ and $
\partial D-\{0\}$ is of class 
${\cal C}^1$.  Let $g:\R^d\rightarrow\R^d$ be a locally Lipschitz vector field, of 
unit length on $\partial D$, such that $g(x)=\bar {g}(x)$ for 
$x\in\partial D\cap\overline {B_{\delta}(0)}$ and, denoting by 
$n(x)$ the unit, inward normal at $x\in\partial D$, it holds  
$\inf_{x\in\partial D-\{0\}}g(x)\cdot n(x)>0$. Then Lemma \ref{th:bdryfunc} 
ensures the existence of a function $F^D\in {\cal C}^2(\overline {
{\cal K}}\cap\overline {B_{\delta}(0)})$ such that 
$\nabla F^D(x)\cdot g\geq c_{F^D}>0$ for every $g\in\bar {G}(x)$, $
x\in\partial {\cal K}\cap\overline {B_{\delta}(0)}$. 
Let $\chi :\R\rightarrow [0,1]$ be a nonincreasing, ${\cal C}^{\infty}$ function such that $
\chi (u)=1$ for 
$u\leq\sup_{x\in\overline {{\cal K}},\,|x|\leq 1}\Phi (x)$ and $\chi 
(u)=0$ for 
$u\geq\inf_{x\in\overline {{\cal K}},\,|x|\geq\delta}\Phi (x)$ . Then the function 
\[F(x):=\chi (\Phi (x))\,F^D(x)+\big(1-\chi (\Phi (x))\big)F^{*}(
x)\]
has the desired properties. 
\hfill $\Box$ \medskip

\end{appendix}
\vskip.5in

\end{document}